\documentclass{amsart}
\usepackage{hyperref}
\usepackage{amssymb, amsmath, amsthm, amsfonts}
\usepackage{verbatim}
\usepackage{bbm}
\usepackage{enumitem}
\usepackage{enumerate}
\usepackage{dsfont}
\usepackage{upgreek}
\usepackage[mathscr]{eucal}
\usepackage{tikz}
\usepackage[normalem]{ulem}
\usepackage{float} 

\usepackage[backgroundcolor=black!5,linecolor=black]{todonotes}
\usepackage{comment}

\newtheorem{thm}{Theorem}[section]

\newtheorem{conj}[thm]{Conjecture}
\newtheorem{prop}[thm]{Proposition}
\newtheorem{cor}[thm]{Corollary}
\newtheorem{lem}[thm]{Lemma}
\numberwithin{equation}{section}
\newtheorem{rem}[thm]{Remark}

\newcommand{\Be}{\begin{equation}}
\newcommand{\Ee}{\end{equation}}
\newcommand\supp{\operatorname{supp}}
 \newcommand{\pq}{(1/p, 1/q)}

\newcommand{\dist}{\operatorname{dist}}

\newcommand{\cE}{\mathcal E}

\begin{document}

\date{\today}
\author[Lee]{Sanghyuk Lee} \address[Lee]{Department of Mathematical Sciences and RIM, Seoul National University, Seoul 08826, Republic of  Korea} \email{shklee@snu.ac.kr}

\author[Roncal]{Luz Roncal} \address[Roncal]{BCAM -- Basque Center for Applied Mathematics, 48009 Bilbao, Spain, and 
Universidad del Pa{\'i}s Vasco / Euskal Herriko Unibertsitatea, 48080 Leioa, Spain, and Ikerbasque, Basque Foundation for Science, 48011 Bilbao, Spain} \email{lroncal@bcamath.org}

\author[Zhang]{Feng Zhang} \address[Zhang]{School of Mathematical Sciences, Xiamen University, Xiamen 361005,
P. R. China} \email{fengzhang@stu.xmu.edu.cn}

\author[Zhao]{Shuijiang Zhao} \address[Zhao]{School of Mathematical Sciences, Zhejiang University, Hangzhou 310058, P. R. China} \email{zhaoshuijiang@zju.edu.cn}

\title[fractal circular maximal  function]{Endpoint estimates for  the fractal   circular maximal  function  and related  local smoothing
}
 
\begin{abstract} 
Sharp $L^p$--$L^q$ estimates for  the  spherical maximal function  over dilation sets of fractal dimensions, including the endpoint estimates, were recently  proved  by Anderson--Hughes--Roos--Seeger. More intricate $L^p$--$L^q$ estimates for the fractal circular maximal function were later established  in the sharp range by  Roos--Seeger, but the endpoint estimates have been left open, particularly when the fractal dimension of the dilation set  lies in $[1/2, 1)$.  In this work, we prove these  missing  endpoint estimates for the circular maximal function. We also study  the closely  related   $L^p$--$L^q$ local smoothing estimates for the wave operator over fractal dilation sets, which were recently investigated by Beltran--Roos--Rutar--Seeger and Wheeler. Making use of a bilinear approach, we also extend the range of $p,q$, for which  the optimal estimate holds. 
\end{abstract}

\maketitle

\section{introduction}
Let $d\ge 2$ and set
\[ A_t f(x)= \int_{\mathbb S^{d-1}} f(x-ty) d\sigma (y), \]
where $\sigma$ is the normalized surface measure on the unit sphere $\mathbb S^{d-1}$.  $L^p$ boundedness of the spherical maximal function $f\mapsto \sup_{t>0}|A_t f|$ is now well understood.  Stein \cite{stein1976} (for $d\ge 3$) and Bourgain \cite{B1} (for $d=2$, where  the maximal operator is referred to as `the circular maximal function') proved that the operator is bounded on  $L^p$ if and only if $p>d/(d-1)$.  Bourgain \cite{B3} also proved a restricted weak type  $(\frac{d}{d-1}, \frac{d}{d-1})$ estimate at the endpoint for $d\ge3$. In contrast, such a restricted weak type estimate fails in two dimensions \cite{STW}, although it does hold for radial functions \cite{Le}.

In this paper, we are concerned with the spherical maximal function  defined by restricted dilation set.  For $E\subset [1,2]$, we consider 
\[ 
M_E f(x)=\sup_{t\in E}  |A_t f(x)|,
\]
and our main focus is on $L^p$--$L^q$ estimates with $p<q$.  When  $E$ consists of a single point (i.e.,  $M_E$ is the average over a single sphere), it is known \cite{littman1973} that the  $L^p$--$L^q$ boundedness  holds if and only if $(1/p,1/q)$ belongs to the closed triangle determined by $(0,0), (1,1)$, and  $(\frac{d}{d+1},\frac{1}{d+1})$. When $E=[1,2]$, $L^p$--$L^q$ estimates for $(1/p,1/q)$  in the interior of  the convex hull of the points 
\[Q_1= (0,0), \ \ Q_{2}=\Big(\frac{d-1}{d},  \frac{d-1}{d}\Big), \ \  Q_{3}= \Big(\frac{d-1}{d}, \frac{1}{d}\Big), \ \ Q_{4}=\Big(\frac{d(d-1)}{d^2+1}, \frac{d-1}{d^2+1}\Big)\] 
  are due to Schlag \cite{S} ($ d=2 $) and Schlag--Sogge \cite{SS} ($ d\geq 3 $).    Note that $Q_2=Q_3$ when $d=2$.  Concerning  the remaining endpoints on the boundary, the first author \cite{L} showed that $M_E$ is of restricted weak type at $Q_4$\footnote{By this, we mean that $M_E$ is of restricted weak type $(p,q)$ for  $(1/p,1/q)=Q_4$}  when $d\ge 2$, and also at $Q_3$ in dimensions $d\ge 3$. 
  In particular, in two dimensions, the restricted weak type at $Q_4$ and  $L^p$–$L^q$ bounds for $(1/p, 1/q)$ in the open line segment $(Q_4, Q_3)$ were established using the bilinear cone restriction estimate due to Wolff \cite{W} (see also \cite{T}).
  
When the dilation set $E$ lies  between the two extreme cases mentioned above,  $L^p$ boundedness of $M_E$  was characterized by Seeger--Wainger--Wright \cite{SWW} in terms of  the upper Minkowski dimension $\dim_M E $ of $E$ (see \eqref{eq:Md} for the definition). Nevertheless, the notion of upper Minkowski dimension alone turned out to be insufficient for determining  the sharp range of $L^p$--$L^q$ boundedness. In  \cite{AHRS} and \cite{RS},   $L^p$--$L^q$  estimates for $M_E$ were studied by using the   Assouad dimension $\dim_A E$   of $E$  (see \eqref{eq:Ma} below) and its variants.   More specifically,  the authors sought to characterize the type set 
\[  T(E)=\{ (1/p, 1/q)\in [0,1]^2: \| M_E\|_{L^p\to L^q} <\infty  \}  \]
in terms of $\dim_M E $, $\dim_A E$, and related notions.

To describe the previous results, we recall some definitions. 
For $\delta > 0$, let $\mathcal N(E,\delta)$ denote  the minimal number of open intervals of length $\delta$ that is required to cover $E$. The (upper) Minkowski dimension $\dim_M E$  is defined by 
\begin{equation}
\label{eq:Md}
\dim_M E = \inf \big\{a > 0 : \exists c > 0 \text{ such that  } \forall \delta \in (0,1),\, \mathcal N(E,\delta) \leq c\delta^{-a}\big\}.
\end{equation}
Let $I$ denote an interval. 
The Assouad dimension $\dim_A E$ (\cite{As})  is defined by
\begin{equation}
\label{eq:Ma}
\begin{aligned}
\dim_A E = \inf \big \{a > 0 : &\, \exists c > 0 \text{ such that }  \\     &\forall I\subset [1,2],  \, \forall \delta \in (0,|I|), \,  \mathcal N(E \cap I, \delta) \leq c\delta^{-a}|I|^a\big\}.
\end{aligned}
\end{equation}

Noting that $0 \leq \dim_M E \leq \dim_A E \leq 1$, we set
\[   \mu=\dim_M E, \quad  \alpha=  \dim_A E.\]   
Following \cite{AHRS} and \cite{RS}, we define  points $Q_1, Q_{2,\mu}, Q_{3,\mu}, Q_{4,\alpha}\in [0,1]^2$  by  
\[
Q_1 = (0,0), \quad Q_{2,\mu} = \left(\frac{d-1}{d-1+\mu},  \frac{d-1}{d-1+\mu}\right),
\]
\[
Q_{3,\mu} = \left(\frac{d-\mu}{d-\mu+1}, \frac{1}{d-\mu+1}\right), \quad Q_{4,\alpha} = \left(\frac{d(d-1)}{d^2+2\alpha-1}, \frac{d-1}{d^2+2\alpha-1}\right).
\]
Moreover, let $\mathcal Q(\mu,\alpha)$ denote the closed convex hull of the points $Q_1$, $Q_{2,\mu}$, $Q_{3,\mu}$, $Q_{4,\alpha}$. 
Define $\mathcal R(\mu,\alpha)$ as the union of the interior of $\mathcal Q(\mu,\alpha)$ and the line segment from $Q_1$ to $Q_{2,\mu}$, excluding   $Q_{2,\mu}$.   Note $Q_{2,\mu}=Q_{3,\mu}$ if  $\mu=1$ and $d=2$. 

Anderson--Hughes--Roos--Seeger \cite{AHRS}   showed that $\mathcal R(\mu, \alpha)\subset  T(E) $  when $d\ge 3$, and when $d=2$ and $\alpha\le 1/2$.  The  more involved case $d=2$ and $\alpha> 1/2$ was later obtained in \cite{RS}. They actually showed a slightly stronger result: The same $L^p$--$L^q$ boundedness remains valid with $\alpha=\dim_A E$ replaced by the quasi-Assouad dimension of $E$.  These results are essentially sharp in that  $L^p$--$L^q$ bound on $M_E$ generally fails unless $(1/p, 1/q)\in \mathcal Q(\mu,\alpha)$. 
This was shown by considering  the quasi-Assouad regular sets (see, for example,  \cite[p.1081]{RS} for the definition).  Recently,  Beltran--Roos--Seeger \cite{BRS} gave an explicit formula for the closure of the $L_{\rm rad}^p\to L^q$ type set for general subsets $E$, where $L_{\rm rad}^p$ denotes the space of  radial $L^p$ functions.  In the radial setting, some endpoint results are also known to hold  \cite{BRS, Zhao}.

\subsection*{Endpoint estimate for the circular maximal function.}  A natural remaining question  is the endpoint $L^p$--$L^q$ boundedness of $M_E$ for $(1/p, 1/q)\in {\mathcal Q}(\mu,\alpha)\setminus \mathcal R(\mu, \alpha)$.

 However, such questions cannot be properly  addressed using only the Minkowski and Assouad dimensions, due to the nature of their definitions. In fact, Anderson et al. \cite{AHRS} obtained endpoint estimates under more quantitative assumptions on the set $E$.
 More precisely, 
the set   $E$ is said to have bounded $\mu$-Minkowski characteristic if 
 \[ \sup_{0<\delta<1}  \delta^\mu \mathcal N(E, \delta)<\infty.\]
 Similarly, we say that  $E$ has bounded $\alpha$-Assouad characteristic if 
\begin{equation} 
\label{assouad}
\mathcal A^{\alpha}(E)=\sup_{0 < \delta < 1} \sup_{I:\delta \leq |I| \leq 1} \Big(\frac{\delta}{|I|} \Big)^{\alpha} \mathcal N(E \cap I, \delta) < \infty,
\end{equation}
 where the inner supremum is taken over  any subinterval $I\subset [1, 2]$ of length $\ge \delta$. 
  
   The following endpoint results were  proved in \cite{AHRS}. 
 
 \begin{thm}[{\cite[Theorem 3]{AHRS}}]  
\label{thm:AHRS} 
Let $ E \subset [1,2] $.
\begin{itemize}[leftmargin=20pt, labelsep=.5em]
\item[$(i)$] Let $ d \ge 2 $ and $ 0 < \mu < 1 $.\footnote{When $ \mu = 1 $  and $ d\geq 3 $, the restricted weak-type estimates at $ Q_{2,1} $ and $ Q_{3,1} $ are covered by \cite[Theorem 1.4]{L}, since $ M_E f \leq M_{[1,2]} f $.}  
If $ E $ has bounded $ \mu $-Minkowski characteristic, then $ M_E $ is of restricted weak type $ (p,q) $ if $(1/p, 1/q)=Q_{2,\mu}$, or $Q_{3,\mu}$.  

\item[$(ii)$] Let $ d \ge 3 $, $ 0 < \alpha \le 1 $, or $ d = 2 $, $ 0 < \alpha < 1/2 $.  
If $ E $ has bounded $ \alpha $-Assouad characteristic, then $ M_E $ is of restricted weak type $ (p,q) $ if $(1/p, 1/q)=Q_{4,\alpha}$.
\end{itemize} 
\end{thm}

We note that the statement in \cite[Theorem 3]{AHRS} concerning sets $ E $ with bounded $ \mu $-Minkowski characteristic appears to be slightly imprecise. Specifically, the restricted weak-type $ (p,q) $ estimates at $ Q_{2,\mu} $ and $ Q_{3,\mu} $ were in fact established for the full range $ 0 < \mu < 1 $ (see \cite[Lemma 2.4]{AHRS}), not only for $ 0 < \mu < 1/2 $.

 The use of bounded Minkowski  and Assouad characteristics in establishing the endpoint estimates  can be rigorously justified. Indeed,  by the argument in the proof of   \cite[Lemma 5.1]{RS},  one can easily see that
 	\begin{equation*}
 	\mathcal N(E \cap I, \delta)^{1/q}(\delta/|I|)^{\frac{d-1}{2}(1-\frac{1}{p}-\frac{1}{q})}\delta^{\frac{d}{q}-\frac{1}{p}}\lesssim \|M_E\|_{L^{p,1}\to L^{q,\infty}}
 	\end{equation*}
	 for any $ I\subset [1,2] $ with $ \delta\leq|I|$.  Note that $\frac{d}{q}-\frac{1}{p}=0$ and $\frac{d-1}{2}(1-\frac{1}{p}-\frac{1}{q})=\frac\alpha q$ if $(1/p, 1/q)=Q_{4,\alpha}$. 
	  Thus,  for the  restricted weak type estimate at $ Q_{4,\alpha}$ to hold, it is necessary  that $E$ has bounded $\alpha$-Assouad characteristic. Similarly, bounded $\mu$-Minkowski characteristic  is necessary for the restricted weak type estimate at  $Q_{3,\mu}$.

 In light of Theorem~\ref{thm:AHRS} and the preceding discussion,  most of endpoint boundedness results have been  established, except at $Q_{4,\alpha}$ when $d=2$ and 
 $\alpha \ge 1/2$.  More  precisely,  the restricted weak type $(p_\alpha,q_\alpha)$ estimate, where
 \[     \Big(\frac1{p_\alpha}, \frac1{q_\alpha}\Big):=  Q_{4,\alpha} = \Big( \frac{2}{2\alpha+3}, \frac1{2\alpha +3}\Big), \] 
  is still unknown.  This problem was explicitly raised in  \cite[Section 2.5]{RS}.   Our first  result is to prove this missing endpoint estimate in the case $d=2$ and $\alpha \ge 1/2$.

 \begin{thm} 
\label{endpoint} Let $d=2$ and $0<\alpha\le1$. Suppose $E\subset[1,2]$ has  bounded $\alpha$-Assouad characteristic. Then $M_E$ is of restricted weak type $(p_\alpha, q_\alpha)$.  
\end{thm} 

Consequently,  under the additional assumption that  $ E $ has  bounded $\mu$-Minkowski characteristic,  this establishes $L^p$--$L^q$ boundedness of $M_E$ for $\pq$ lying in $(Q_1, Q_{4, \alpha}) \cup ( Q_{4, \alpha}, Q_{3, \mu}) $, where $(Q_1, Q_{4, \alpha})$ and $( Q_{4, \alpha}, Q_{3, \mu}) $ denote the open line segments connecting 
$Q_1$ to $Q_{4, \alpha}$ and $Q_{4, \alpha} $ to $Q_{3, \mu}$, respectively. 

\subsection*{Local smoothing estimate over fractal time set} 

It is well known (\cite{MSS, S, SS, TV2, L}) that the study  of the circular and spherical maximal functions is closely related to smoothing properties of the wave propagator 
\[  e^{it\sqrt{-\Delta}} f(x)= \frac1{(2\pi)^d}   \int e^{i( x\cdot \xi+  t|\xi|)}  \widehat f(\xi) d\xi.\]
(See, in particular, \cite[Remark 1.4]{S} for a discussion on the connection with the circular maximal function.)
The recent developments in \cite{AHRS, RS} related to the fractal spherical maximal function are based on local smoothing estimates  that incorporate  the fractal structure of the underlying set. 

We say that an  arbitrary subset $\mathcal E\subset [1,2]$ is $\delta$-separated if $|t-s|\ge \delta$ whenever $t,s\in \cE$ and $s\neq t$. 
By $E(\delta)$  we denote  a   $\delta$-discretization of $E$, i.e., a maximally  $\delta$-separated subset of $E$.
 Let $\beta\in  C_c^\infty ((1/2, 2))$  be such that $\sum_{j=-\infty}^\infty \beta(2^{-j} t)=1$ for $t>0$. 
 Let $P_j$ be the standard Littlewood--Paley decomposition given by 
\[ \widehat{P_j f}(\xi)=\beta(|\xi|/2^j)  \widehat f(\xi), \quad j\ge 1.  \]  

For $1\le p\le q\le \infty$, assuming  that the set $E$ has bounded $\alpha$-Assouad characteristic, i.e.,   \eqref{assouad} holds,  we consider  the following form of smoothing estimates, which have played an important role in the previous works:
\Be
\label{smooth-pq}
  \Big( \sum_{t\in  E(2^{-j})} \|  e^{it\sqrt{-\Delta}} {P_j f}\|_{L_x^q(\mathbb R^d)}^q \Big)^\frac1q \le C 2^{sj} 
     \|f\|_p 
  \Ee
  with suitable $s$ depending on $p, q,$ and  $\alpha$.  A precise description of 
$s$  is given below.    A different formulation was presented in \cite{BRRS} by using the  Legendre--Assouad function (see also \cite{Wh}).

   For $j=1,2,3$, let us define $s_j(p,q)$ by 
  \[ 
  \begin{aligned}
   s_1(p,q) &\textstyle=\frac{d-1}{2}+\frac{1}{p}-\frac{d}{q}, 
   \\
   s_2(p,q) &\textstyle = \frac{d+1}{2}\big(\frac{1}{p}-\frac{1}{q}\big)+\frac{{\alpha}}{q},  
    \\
 s_3(p,q) &\textstyle =\frac{d}{p}-\frac{1-{\alpha}}{q}-\frac{d-1}{2} .
   \end{aligned} 
   \]
We  set 
\[
s_c(p,q)=\max  \big\{s_1(p,q),\,s_2(p,q),\, s_3(p,q)\big\}.
\]
  
  It is not difficult to see  that 
   \eqref{smooth-pq} holds only if   $s\ge s_c(p,q)$ (see Proposition \ref{prop:neces}  below). In fact,  one can show  that there is a set $E$ of bounded $\alpha$-Assouad characteristic, for which 
     \eqref{smooth-pq}  fails  if $s< s_c(p,q)$. 
   In particular, when $p=2$, the  estimate \eqref{smooth-pq} with the optimal exponent  $s=s_c(p,q)$ was  established in \cite{Wh} and \cite{BRRS} except for the endpoint case $q=2(d-1+2{\alpha})/(d-1)$ 
   (see also \cite{AHRS} for an earlier related result).

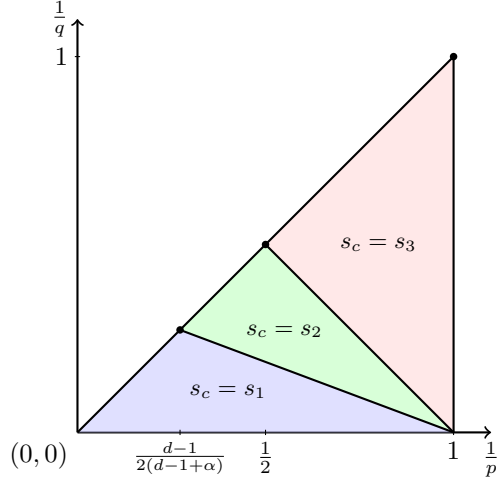
\begin{figure}[t] 
	\centering 
	\begin{tikzpicture}[scale=0.5]
		\draw[thick,->]  (0,0) -- (11,0);
		\draw[thick,->] (0,0) -- (0,11);
		\node[below] at (11,0) {$\frac{1}{p}$};
		\node[left] at (0,11) {$\frac{1}{q}$};
		
		\fill[blue!20,opacity=0.6] (0,0) -- (10,0) -- (30/11,30/11) -- cycle;
		\node[font=\small] at (4,1.1) {$s_c=s_1$};
		
		\fill[green!25,opacity=0.6] (30/11,30/11) -- (10,0) -- (5,5) -- cycle;
		\node[font=\small] at (5.5,2.7) {$s_c=s_2$};
		
		\fill[red!15,opacity=0.6] (5,5) -- (10,0) -- (10,10) -- cycle;
		\node[font=\small] at (8,5) {$s_c=s_3$};
		
		\draw[thick] (0,0) -- (10,10);        
		\draw[thick] (5,5) -- (10,0);         
		\draw[thick] (10,0) -- (10,10);       
		\draw[thick] (30/11,30/11) -- (10,0); 
	    
	    \node[below left] at (0,0) {$(0,0)$};    
	    
	    \draw (10,-2pt) -- (10,2pt);             
	    \node[below] at (10,0) {$1$};
	  
	    \draw (-2pt,10) -- (2pt,10);             

	    \node[left] at (0,10) {$1$};
	    \draw (5,-2pt) -- (5,2pt);               
	    \node[below] at (5,0) {$\frac{1}{2}$};
	    \node[below,font=\footnotesize] at (30/11,0) { $\frac{d-1}{2(d-1+\alpha)}$};
	    
	    \node[circle, fill=black, inner sep=1.0 pt] at (5,5) {};                                
	    \node[circle, fill=black, inner sep=1.0 pt] at (10,10) {};                                
	    \node[circle, fill=black, inner sep=1.0 pt] at (30/11,30/11) {};
	    \draw (30/11,-2pt) -- (30/11,2pt);
	\end{tikzpicture}
	\caption{ Sharp regularity exponent $s_c(p,q)$ for the fractal local smoothing estimate \eqref{smooth-pq} when $d=2$ and $\alpha=5/6$.} 
	\label{Fig1} 
\end{figure}

For $t\in [1,2]$, the fixed-time estimate
\Be 
\label{single}   \| e^{it\sqrt{-\Delta}} P_jf\|_{L_x^q(\mathbb R^d)}\lesssim 2^{j (s_1(p,q)\vee s_1(q',p'))} \|f\|_p  
\Ee
follows from the estimates for the kernel of the operator  $f\mapsto e^{it\sqrt{-\Delta}} P_jf$ and duality. Here $a\vee b=\max\{a, b\}$.  
Moreover, the estimate \eqref{single} is sharp in that the regularity exponent $s_1(p,q)\vee s_1(q',p')$ cannot be replaced by any smaller one. 
In particular,  the sharp $L^p_s$--$L^p$ estimate for the wave operator  $e^{it\sqrt{-\Delta}}$ with  $1<p<\infty$ and $s=s_1(p,p)\vee s_1(p',p')$ goes back to  Miyachi \cite{Mi}  and Peral \cite{Pe}.  Note that 
\[
s_c(p,q)=
\begin{cases}
s_1(p,q), \qquad   (1-\frac1p)\ge \frac{d-1+2{\alpha}}{(d-1)q}, 
\\
s_2(p,q), \qquad  \frac1q\le (1-\frac1p)\le \frac{d-1+2{\alpha}}{(d-1)q}, 
\\
s_3(p,q), \qquad (1-\frac1p) \le \frac1q. 
\end{cases} 
\]
(See  Figure \ref{Fig1}.)
  As will be seen later, the estimate \eqref{smooth-pq} is not of particular interest when $s=s_c(p,q)=s_3(p,q)$ since in this case it follows from the fixed time estimate \eqref{single}. 
  
The estimate \eqref{smooth-pq} becomes most interesting  when $s=s_c(p,q)=s_1(p,q)$. In this case,  \eqref{smooth-pq} represents  a substantial improvement over the fixed-time estimate \eqref{single} in the sense that  the same sharp bound of  \eqref{single}  continues to control the left hand side of  \eqref{smooth-pq}  even if there are many  $t\in E(2^{-j})$. 
Moreover, the most significant case occurs  when $s_1(p,q)=s_2(p,q)$, that is to say, 
\Be 
\label{s-line}
 (d-1)\Big(1-\frac1p\Big)=\frac{d-1+2{\alpha}}q.
\Ee
Once one obtains the estimate    \eqref{smooth-pq}  for the case $s_1(p,q)=s_2(p,q)$, then the other estimates basically  follow by interpolation.  Consequently, we are led to conjecture the following. 

\begin{conj}  Let $0<\alpha\le 1$, and $p,q$ satisfy \eqref{s-line}. Suppose that  $ E\subset [1,2] $ has  bounded $\alpha$-Assouad characteristic. Then, the estimate \eqref{smooth-pq} with $s=s_c(p,q)$ 
holds for $2(d-1+\alpha)/(d-1)< q< \infty$. 
\end{conj}

For $p=q=2(d-1+\alpha)/(d-1)$,  failure of the estimate \eqref{smooth-pq}  is not difficult to show  (see Proposition \ref{prop:neces} below).  In this case,  the estimate does not hold without an additional $2^{\epsilon j}$ loss in the bound.

A related conjecture concerning  the  $L^p$  estimate for radial functions  was recently verified by Beltran--Roos--Rutar--Seeger  \cite[Theorem 1.1]{BRRS}.  For nonradial general function, the conjecture was only verified for $ 2(d-1+2\alpha)/(d-1)< q$ in \cite{BRRS} and \cite{Wh}, independently.  The approach in the previous works is essentially based on a $TT^*$ argument.

We  extend the range of sharp boundedness by employing bilinear estimates. Since we aim to control $f\mapsto e^{it\sqrt{-\Delta}} P_jf$ over a restricted set, it is natural to consider 
bilinear estimates  which exploit the sparseness of the $t$-set $E$ (see Theorem \ref{bilinear-stronger}). It is very likely that the range can be further extended by making use of multilinear or the refined arguments currently available. This will be pursued elsewhere.

\begin{thm} 
\label{improved} 
Let $0<\alpha\le 1$, and $p,q$ satisfy \eqref{s-line}. Suppose that  $ E\subset [1,2] $ has  bounded $\alpha$-Assouad characteristic. Then, the estimate \eqref{smooth-pq} with $s=s_c(p,q)$ 
holds for  $q>\frac{2(d-1+2\alpha)^2-4\alpha^2}{(d-1)(d-1+2\alpha)}$.
\end{thm} 

In particular, Theorem \ref{improved} gives the  estimate \eqref{smooth-pq} with the optimal exponent  $s=s_c(p,q)$  for  $q=2(d-1+2\alpha)/(d-1)$,  i.e.,  $p=2$.  For this case,  the  weak type  $L^{2}\to L^{q,\infty}(\ell^2)$ estimate  \cite[Proposition 6.3 $(ii)$]{BRRS}  and the strong type  $L^{2}\to L^{q}(\ell^2)$ estimate with $\epsilon$-loss in regularity  were respectively obtained in 
\cite{BRRS} and \cite{Wh}. If a $2^{\epsilon j}$-loss is allowed, the range can be further extended by making use of the known sharp local smoothing estimate (Section \ref{sec:sss}). Moreover, when $d\ge 4$, the range of $q$ can be enlarged thanks to the optimal local smoothing estimate due to  Heo--Nazarov--Seeger \cite{HNS} (see Remark \ref{hns}). 

\subsubsection* {Organization}  We first prove \eqref{smooth-pq} for an extended range $p,q$ by directly applying the bilinear restriction estimates for the cone, which are due to Wolff \cite{W} and Tao \cite{T}. 
The consequent estimates are enough to show the endpoint estimate in Theorem \ref{endpoint}. This is done in Section \ref{sec:localsmoothing}.    
Subsequently, in Section \ref{sec:improved},  we further extend the $p,q$ range  of the estimate \eqref{smooth-pq} with $s=s_c(p,q)$.  For the purpose,  making use of the induction on scales argument, we  obtain a bilinear restriction estimate over sparse time sets  to prove Theorem \ref{improved}.  In Section \ref{sec:lower}, we discuss  the 
sharpness of the regularity exponent $s_c(p,q)$. 

\section{Local smoothing estimates over discrete temporal sets}
\label{sec:localsmoothing}

Let $\cE\subset [1,2]$   be a finite  discrete set. 
For $F$ defined on $\mathbb R^d\times  \cE$ and $q\in [1, \infty]$, we denote 
\[  \|F\|_{L^q( \mathbb R^d\times\cE)}=\Big( \sum_{t\in  \cE} \| F(x,t)\|_{L_x^q(\mathbb R^d)}^q \Big)^\frac1q  \] 
if $q\neq \infty$ and $ \|F\|_{L^\infty( \mathbb R^d\times \cE)}= \max_{t\in  \cE} \| F(\cdot,t)\|_{L_x^\infty(\mathbb R^d)} $.  As in \cite[Section 2.3]{RS}, 
 define
	\Be
\label{agamma}
\mathcal A_{\alpha}(\mathcal E;  \delta):=\sup_{I: \delta\le  |I|  \le 1}\Big(\frac{\delta}{|I|}\Big)^{\alpha}\mathcal N(\mathcal{E}\cap I,\delta),
\Ee
where $I$ denotes subintervals of $[1,2]$. 
Note that $\mathcal N(\mathcal{E}\cap I,\delta)=\# (\mathcal E\cap I)$ if $\mathcal E$ is a $\delta$-separated set. 

 The following local smoothing estimate will be used in the proof of Theorem \ref{endpoint}.

\begin{prop}\label{p-ls} Let $(d-1)/(d+1)< \alpha \le 1$ and $N\ge 1$.  Let 
\[ 
\mathbb A_N=\{\xi: N/2\le |\xi|\le 2N\}.
\] Let $ \cE\subset [1,2]$ be $1/N$-separated. 
 Suppose that $\widehat f$ is supported in $\mathbb A_N$. 
Then, for $q> q_\ast:= \frac{2(d^2+(\alpha+1) d+\alpha-2)}{d^2-1}$ and $p$ given by \eqref{s-line},   we have the estimate 
\Be 
\label{localsmoothing-2d1} 
\| e^{\pm it\sqrt{-\Delta}} f\|_{L^q( \mathbb R^d\times  \cE)} \le C \mathcal A_{\alpha}^\frac1q N^{\frac{d-1}2+\frac1p-\frac{d}{q}} \|f\|_p,
\Ee  
where $
\mathcal A_{\alpha}=\mathcal A_{\alpha}(\cE; 1/N)$. 
\end{prop}

When $q=\infty$, the estimate coincides with the fixed time estimate \eqref{single}.  Thanks to the fixed-time estimate \eqref{single}, the estimate \eqref{localsmoothing-2d1} trivially  holds  when $N$ is bounded by a constant. Thus, we always assume that $N$ is sufficiently large.  Since  $p$, $q$ satisfy \eqref{s-line}, the exponent of $N$ can be written in a more concise way, that is to say, 
\[ \frac{d-1}2+\frac1p-\frac dq  = \frac{d+1}2-\frac{d^2-1+ 2\alpha}{(d-1)q} .\]
However, any  expression (of the exponent) of the form  $c+\frac bp+ \frac{a}{q}$, which remains invariant\,\footnote{Indeed, let $B(p,q)=\mathcal A_{\alpha}^\frac1q N^{\frac{d-1}2+\frac1p-\frac{d}{q}}$. If we have  the estimate \eqref{localsmoothing-2d1} for $(p,q)=(p_j, q_j)$ for $j=0,1$, via interpolation 
we get  \eqref{localsmoothing-2d1}  with bound $B(p_0,q_0)^{1-\theta}B(p_1,q_1)^\theta=B(p,q)$ with  $1/p=(1-\theta)/p_0+ \theta/p_1$ and $1/q=(1-\theta)/q_0+ \theta/q_1$.}   while  we interpolate  the estimate \eqref{localsmoothing-2d1}, is acceptable for our purpose.  So, we do not insist on 
writing the exponent in a unique  way.

We also make use of the next lemma, which is a multilinear extension of Bourgain's summation trick  \cite{B3, CSW}. 
\begin{lem}[{\cite[Lemma 2.6]{L}}]
\label{B-Sum-Tri} Let $n\ge 1$. 
Let  $1\leq p_0^k,p_1^k\le \infty$, $k=1,\cdots, n,$ and  $1\leq q_0, q_1\le\infty$. Suppose that $\{T_j\}_{j=-\infty}^\infty$ is a sequence of $n$-linear (or sublinear) operators such that
\begin{equation*}
\|T_j(f^1,\cdots,f^n)\|_{L^{q_\ell}}\leq M_\ell  2^{j (-1)^\ell \varepsilon_\ell }\prod_{k=1}^n\|f^k\|_{L^{p_\ell^k}}, \quad \ell=0, 1
\end{equation*}
for some  $\varepsilon_0$, $\varepsilon_1>0$.  Let $\theta=\varepsilon_1/(\varepsilon_0+\varepsilon_1)$, $1/q=\theta/q_0+(1-\theta)/q_1$, and $1/p^k=\theta/p_0^k+(1-\theta)/p_1^k$.
Then, we have 
\begin{equation*}
\|(\sum T_j)(f^1,\cdots,f^n)\|_{L^{q,\infty}}\leq CM_1^{\theta}M_2^{1-\theta}\prod_{k=1}^{n}\|f^k\|_{L^{p^k,1}}.
\end{equation*}
\end{lem}

  \subsection{Proof of Theorem \ref{endpoint}}
Since the case $\alpha<1/2$  was already obtained in \cite{AHRS}, we may assume $\alpha\ge 1/2$. 

Let us  set $P_0 =1- \sum_{j\ge 1} P_j f$.  For $j\ge 0$,  we also set
\[  M_j f(x)=\sup_{t\in E} |A_t P_jf(x)|.\]
Then, it is easy to show $\|M_0 f\|_q\le C\|f\|_p$ for $1\le p\le q\le \infty$. Thus, instead of $M_E$,  it is sufficient to show that  
\[  f\mapsto    {\widetilde Mf}:= \sum_{j=1}^\infty  M_j f\] 
is of restricted weak type $(p_\alpha, q_\alpha)$. 
We now claim that
 \Be 
\label{max} 
\| M_j f \|_{L^q( \mathbb R^2)} \le C    2^{(1-\frac{2\alpha+3}q)j} \|f\|_p
\Ee   
provided that $(1-\frac1p)=\frac{1+2\alpha}q$  and $q> (8+6\alpha)/3$. Note that  $(8+6\alpha)/3< q_\alpha= 2\alpha+3$. 
Therefore, once we have the estimate \eqref{max}, we get the restricted weak type $(p_\alpha, q_\alpha)$ estimate using  Lemma~\ref{B-Sum-Tri} with $n=1$.

We now proceed to prove the claimed estimate \eqref{max}. To this end, we require some preliminary observations. Since we are working with $P_jf$ whose Fourier transform is supported in $\mathbb A_{2^j}$, the function  $t\mapsto e^{it\sqrt{-\Delta}} P_jf(x)$ has the \textit{locally constant property at 
scale $2^{-j}$}, that is to say,  $t\mapsto e^{it\sqrt{-\Delta}} P_jf(x)$ behaves as if it were a constant within any interval of length  $2^{-j}$ (see \eqref{local-sim} below).  To make it  precise, for a function $F$,  we define a scattered translation sum of $F$ as follows:  
\Be 
\label{def-sj} 
\mathfrak S F(x) =  \sum_{\mathbf k\in \mathbb Z^d}   c_{\mathbf k}  |F(x+ \tau_{\mathbf k})|,   
\Ee
where $c_{\mathbf k}=C(1+|\mathbf k|)^{-100d}$ for a  constant $C$ and $\mathbf k\mapsto\tau_{\mathbf k}$  is a function from $\mathbb Z^d$ to $\mathbb R^d$. 
For our purpose, specifying the function $\mathbf k\mapsto\tau_{\mathbf k}$ is not necessary so that it may vary in different occurrences of the operator $\mathfrak S$.  For example, it may depend on $j$ when $ F(x)=e^{it\sqrt{-\Delta}} P_jf(x)$.   Nevertheless, we shall abuse notation and continue to denote this operator by 
$\mathfrak S$  as this ambiguity is harmless in our context.

\begin{rem}
\label{rem:s-sum}
Let $B, B'$ be Banach spaces and  $T$ be a translation invariant sub-linear operator from $B$ to $B'$  such that $\|T\|_{B\to B'}\le \mathcal B$. 
Suppose that the norm on  $B$ or $B'$ is  translation invariant. Then, 
by the rapid decay of $c_{\mathbf k}$ it is clear that there is a constant $C$ such that $\|\mathfrak ST f\|_{B'}\le C\mathcal B \|f\|_{ B}$.  
\end{rem}

Using the operator $\mathfrak S$, we obtain the  following lemma,  which allows us to treat   $|e^{it\sqrt{-\Delta}} P_jf(x)|$ and  $|e^{it'\sqrt{-\Delta}} P_jf(x)|$ as if they satisfy  
\Be
\label{local-sim}
 |e^{it\sqrt{-\Delta}} P_jf(x)|\sim |e^{it'\sqrt{-\Delta}} P_jf(x)|
 \Ee 
 when $|t-t'|\le 2^{-j}$ (see \cite{L3} for a similar idea used in a different setting).

\begin{lem}[{Locally  constant property}]
\label{lcp} 
Let $\,t, t'\in [1,2]$ and $\supp \widehat f\subset \mathbb A_{2^j}$. Suppose $|t-t'|\le 2^{-j}$. Then
\[    |e^{\pm it\sqrt{-\Delta}} f(x)| \le   \mathfrak S(e^{\pm it'\sqrt{-\Delta}} f)( x) \le  \mathfrak S^2(e^{\pm it\sqrt{-\Delta}} f)( x) . \]  
\end{lem}

\begin{proof} It is sufficient to show the first inequality. Changing variables $\xi\to 2^j\xi$, we have 
\[   e^{\pm it\sqrt{-\Delta}} f(x)=  
C_d  2^{dj}   \int e^{i2^j( x\cdot \xi\pm  t'|\xi|)} m(\xi) \widehat f(2^j \xi) d\xi, \]  
where  $C_d=(2\pi)^{-d}$ and 
 \[ m(\xi) = \beta( \xi)  e^{\pm i2^j(t-t')|\xi|}.\]  
 Recalling that $\beta\in  C_c^\infty((1/2, 2))$,  we expand $m$ in Fourier series over the set $[-\pi, \pi]^d$ to have $m(\xi)=  \sum_{\mathbf k\in  \mathbb Z^d}  d_{\mathbf k, \pm} (t, t')  
e^{i \mathbf k\cdot \xi}.  $ Since $|\partial_\xi^\alpha  m|\le C_\alpha$ for any $\alpha$,  for any $M>0$ there exists a constant $C$,  independent of $t,t'$,  such that $|d_{\mathbf k, \pm} (t, t')|\le C(1+|\mathbf k|)^{-M}$ with $C$.   
Thus,  putting the expansion of $m$ in the integral and reversing the change of variables $\xi\to 2^j\xi$, we get the desired inequality. 
\end{proof}

Thanks to  the locally constant property at 
scale $2^{-j}$, 
it suffices to consider a  maximally $2^{-j}$-separated subset of $E$. Let $E_j$ be such a  set.  Then,  for $j\ge 1$  it is clear that $ \# (E_j \cap I)\le \mathcal N(E\cap I, 2^{-j})$ 
for any interval $I\subset [1,2]$ with $|I|\ge 2^{-j}$.  Thus, recalling \eqref{agamma} and \eqref{assouad},  
 we have
\Be 
\label{discrete} 
\mathcal A_{\alpha}(E_j; 2^{-j})=\sup_{ I\subset [1,2],\,  |I|\ge 2^{-j}} \Big(\frac{2^{-j}}{|I|} \Big)^{\alpha} \#(E_j \cap I) \le  \mathcal A^{\alpha}(E)
\Ee
for $j\ge 1$ and a positive constant $C$.

To show the estimate \eqref{max}, we now use  the asymptotic expansion of $\widehat{\sigma}$ to  have 
\[ A_t P_jf(x)= \sum_{\pm} c_\pm  2^{-j/2} e^{\pm it\sqrt{-\Delta}} P_j \tilde f 
\]
where $\|\tilde f\|_p\le C\|f\|_p$.  Since every point in $E$ is within a $2^{-j}$ neighborhood of some point of $E_j$, 
using Lemma \ref{lcp},  we see that 
\begin{align*}
 M_j f(x)= \sup_{t\in E} |A_t P_jf(x)| & \lesssim   2^{-j/2}   \sum_{\pm}  \max_{t\in E_j}  \mathfrak S(e^{\pm it\sqrt{-\Delta}} P_j \tilde f\,)(x)
 \\
 & \lesssim   2^{-j/2}  \sum_{\pm}  \Big(\sum_{t\in E_j}   \big(\mathfrak S(e^{\pm it\sqrt{-\Delta}} P_j \tilde f\,)(x)\big)^q\Big)^{1/q} . 
 \end{align*} 
Thus, we have 
\[ \|M_j f\|_{L^q( \mathbb R^2)}\lesssim  2^{-j/2}  \sum_{\pm}  \| \mathfrak S(e^{\pm it\sqrt{-\Delta}} P_j \tilde f\,)\|_{L^q( \mathbb R^2\times E_j)} .\]

Now recall \eqref{def-sj}.   By the triangle inequality, the estimate \eqref{max} follows (see Remark~\ref{rem:s-sum}) if we obtain  
\[\| e^{\pm it\sqrt{-\Delta}} f\|_{L^q( \mathbb R^2\times E_j)} \le C 2^{(\frac{3}2-\frac {3+2\alpha}q )j} \|f\|_p\]
when  $(1-\frac1p)=\frac{1+2\alpha}q$  and $q> (8+6\alpha)/3$. Since $E$ has bounded $\alpha$-Assouad characteristic,    \eqref{discrete}  implies   $\mathcal A_{\alpha}(E_j; 2^{-j}) \le C$ for a constant $C$. Therefore, the desired estimate  is an immediate consequence of  Proposition \ref{p-ls} with $d=2$. \qed

\subsection{Proof of Proposition \ref{p-ls}} 

To prove Proposition \ref{p-ls}, we follow the strategy employed  in \cite{L} (see also \cite{CLW}), which  relies on  the bilinear restriction estimate for the cone. 
Let 
\Be \label{rest-o} \mathcal R^\ast g(x,t)= \int  e^{i( x\cdot \xi+  t|\xi|)}  g(\xi) d\xi.\Ee
We now recall the following result, due to Wolff \cite{W} and Tao \cite{T}.   

\begin{thm}  
\label{thm:cone}
Let $\Theta_1, \Theta_2$ be subsets of $\mathbb S^{d-1}$  such that $\dist(\Theta_1, \Theta_2)\ge  c$ for a constant $c>0$. Suppose  that 
\begin{align}
 \supp g_j \subset   \{ \xi\in \mathbb A_1:    \xi/|\xi|\in \Theta_j \},  \quad j=1,2. 
  \end{align}
    Then, for $q\ge q_\circ=q_\circ(d):=2(d+3)/(d+1)$, we have 
\Be 
\label{bilinear0}
\| \prod_{j=1}^2 \mathcal R^\ast g_j \|_{L^{\frac q2}( \mathbb R^{d+1})} \le C 
\prod_{j=1}^2  \|g_j\|_{L^2}. 
\Ee
\end{thm}

 The estimate   \eqref{bilinear0}  for  $q>q_\circ$  is due to  Wolff \cite{W} and the endpoint case $q=q_\circ$ was obtained  by Tao \cite{T}.

As a consequence of Theorem \ref{thm:cone}, we obtain the following proposition. 

\begin{prop}
\label{thm:bi}  Let   $1/\sqrt N \le \theta \le 1/10$.   Let $\Theta_1, \Theta_2$ be subsets of $\mathbb S^{d-1}$ of diameter $\theta$ such that $\dist(\Theta_1, \Theta_2)\sim  \theta$. Suppose  that 
\begin{align}
\label{supp-con}
 \supp \widehat f_j\subset \Lambda_j &:=\{ \xi:  \xi/|\xi|\in \Theta_j, \ N/2\le  |\xi|\le 2N \},   \quad j=1,2. 
   \end{align}
  Then, for $p, q$ satisfying  $q_\circ/q=2/p$ and $q_\circ\le q\le\infty$, we have 
\Be 
\label{bilinear}
\| \prod_{j=1}^2 e^{it\sqrt{-\Delta}} f_j \|_{L^{\frac q2}( \mathbb R^d\times  I_\circ)} \lesssim  
N^{2(\frac{d-1}2+\frac{1}{p}- \frac{d+1}q)}
 \theta^{ 4(\frac{d-1}2-\frac{d-1}{2p}- \frac{d+1}{2q})} \prod_{j=1}^2 \|f_j\|_{L^p}.
\Ee
Here $I_\circ=[1,2]$.  Furthermore, if $\theta\sim 1/\sqrt N$, then \eqref{bilinear} holds without the assumption that $\dist(\Theta_1, \Theta_2)\sim \theta$. 
\end{prop}

The last statement in the theorem  can be easily shown by using Bernstein's inequality.  When $\theta\sim 1$, $N\sim 1$, and $p=2$,  the estimate \eqref{bilinear} follows from \eqref{bilinear0} by Plancherel's theorem.  The estimates  for the case $1/\sqrt N \ll \theta \le 1/10$ can be obtained  by rescaling (see, for instance, \cite[Lemma 2.5]{L}).

\begin{proof}[Proof of Proposition \ref{thm:bi}]  
 By rotation we may assume  that the  sets $\Theta_1$ and $\Theta_2$ are contained in a $C\theta$ neighborhood of $e_1$. More precisely, we may assume 
\Be\label{supp0}
 \supp \widehat f_j \subset \Lambda=\{ \xi:   
  \xi_1\in [N/2, 2N], \  |\xi'| \le  C \theta N    \},\quad j=1,2,    
 \Ee
where $\xi=(\xi_1, \xi')\in\mathbb R\times \mathbb R^{d-1}$. 

To prove Proposition \ref{thm:bi}, by interpolation it is sufficient to show \eqref{bilinear} for $(p,q)=(\infty, \infty)$ and 
$(p,q)=(2,q_\circ)$.  We first consider the case  $(p,q)=(\infty, \infty)$, which is easy. Let $\beta_0\in C_c^\infty ((-2^2, 2^2))$ be such that  $\beta_0=1$ on $[-2, 2]$,  and let $\beta_1 \in C_c^{\infty} ((2^{-2}, 2^2))$ be such that $\beta_1=1$ on $[2^{-1}, 2]$.   
Let us set 
\[  \chi^{N, \theta} (\xi)= \beta_1\Big(\frac{\xi_1}{N}\Big)\beta_0\Big(\frac {|\xi'|}{C\theta N}\Big)  \] 
and
\[K_t^{N,\theta}(x)=(2 \pi)^{-d} \int  e^{i (x\cdot \xi+t|\xi|)} \chi^{N,\theta} (\xi) d\xi. \] 
Note $e^{it\sqrt{-\Delta}} f_j=K_t^{N,\theta}\ast f_j$ for $j=1,2$. It is not difficult to see that   $\|K_t^{N,\theta}\|_{L^1}\le C$  if $\theta\le N^{-1/2}$. Thus, dividing the support of  $\chi^{N,\theta}$ into as many as $(\theta N^\frac12)^{d-1}$   sectors of  angular diameter $N^{-1/2}$, we have $\|K_t^{N,\theta}\|_1\le C(\theta N^\frac12)^{d-1}$ for $t\in I_\circ=[1,2]$. Hence 
$\|e^{it\sqrt{-\Delta}} f_j\|_\infty \le C(\theta N^\frac12)^{d-1}\| f_j\|_\infty$. This gives \eqref{bilinear} for $(p,q)=(\infty, \infty)$.

To show 
\eqref{bilinear} for $(p,q)=(2, q_\circ)$, we use the bilinear endpoint restriction estimate for the cone  (Theorem \ref{thm:cone}). Since $\dist(\Theta_1, \Theta_2)\sim  \theta$, 
by decomposing  $\supp \widehat f_2$ in finite angular sectors and additional rotation, we may assume  that
\Be
\label{sep}
\begin{aligned}
	\supp \widehat f_1 &\subset \{ \xi\in \Lambda:  |\xi'| \sim \theta N       \},    
\\ 
	\supp \widehat f_2  &\subset \{ \xi\in \Lambda:  |\xi'|\le c \theta N    \}
\end{aligned}
\Ee 
for a small constant $c>0$.
 By  scaling  (i.e., $\widehat f_j  \to   \widehat f_j (N\cdot)$), we may assume 
$N=1$. Consequently, it is enough to show 
\Be \label{reduced} \| e^{it\sqrt{-\Delta}} f_1  e^{it\sqrt{-\Delta}} f_2 \|_{L^{r}( \mathbb R^d\times \mathbb R)} \lesssim    \theta^{d-1-\frac{d+1}r} \|f_1\|_{L^2} \|f_2\|_{L^2}\Ee
with $r:=q_\circ/2$  while \eqref{sep}  holds with $N=1$.

 Let us set
\Be \label{extension}
\mathcal T h(x,t)= \iint e^{i(x_1\rho+ x'\cdot \eta+ t  |\eta|^2/2\rho)}  h(\rho, \eta) d\eta d\rho. \Ee
By changing variables 
\[
(x, t) \mapsto \big( 2^{-1/2}(x_1 - t),\, x',\, 2^{-1/2}(x_1 + t) \big), 
\]
  the phase function $(x,t)\cdot (\xi, |\xi|)$ for $e^{it\sqrt{-\Delta}} f$ is transformed to $(x,t)\cdot(2^{-1/2}(\xi_1+|\xi|), \xi', 2^{-1/2}(|\xi|-\xi_1))$. By an additional change of variables $(2^{-1/2}(\xi_1+|\xi|), \xi')\to (\rho, \eta)$,  we note that 
\[ 2^{-1/2}(|\xi|-\xi_1)=|\eta|^2/2\rho,\] thus by  Plancherel's theorem and discarding harmless smooth factor resulting from the change of variables,  the estimate \eqref{reduced} is now equivalent to the estimate 
\Be 
\label{theta} 
\| \mathcal T h_1  \mathcal T h_2 \|_{L^{r}( \mathbb R^d\times \mathbb R)} \lesssim    \theta^{d-1-\frac{d+1}r} \|h_1\|_{L^2} \|h_2\|_{L^2},
\Ee
while 
\Be
\label{sep-scaled}
\begin{aligned}
  \supp h_1 \subset  \mathfrak A_\theta&:= \{ ({\rho,\eta}) : \rho\sim 1,   |\eta| {\sim \theta}    \}, 
\\
     \supp h_2\subset \mathfrak A_\theta'&:= \{ ({\rho,\eta}) : \rho\sim 1,  |\eta|  \le c \theta  \},
\end{aligned}  
\Ee
where $c>0$ is a small constant. 
Clearly, \eqref{theta} with $\theta\sim 1$ is equivalent to the bilinear restriction estimate  \eqref{bilinear0} with $q=q_\circ$.

To show \eqref{theta}  for small $\theta$, we now set
\[\tilde h_j(\rho, \eta)= \theta^{\frac{d-1}2}h_j(\rho, \theta \eta), \quad j=1,2.\]  
Changing variables $\eta\to \theta \eta$,  we have 
\[ (\mathcal T h_1  \mathcal T h_2)(x,t)=  \theta^{d-1}(\mathcal T \tilde h_1  \mathcal T \tilde h_2)(x_1, \theta x',\theta^2 t).   \]
Consequently, 
\[  \| \mathcal T h_1  \mathcal T h_2 \|_{L^{r}( \mathbb R^d\times \mathbb R)} = \theta^{d-1-\frac{d+1}r} \| \mathcal T \tilde h_1  \mathcal T \tilde h_2 \|_{L^{r}( \mathbb R^d\times \mathbb R)}.\]
Since $\supp  \tilde h_1\subset \mathfrak A_1$ and $\supp  \tilde h_2\subset \mathfrak A_1'$, we may apply 
the estimate \eqref{theta} with $\theta=1$ to the right hand side. Since $\|\tilde h_1\|_2 =  \|h_1\|_2$ and $\|\tilde h_2\|_2 =  \|h_2\|_2$, the desired estimate 
 \eqref{theta} follows.

When  $\theta\sim 1/\sqrt N$,  the separation condition $\dist(\Theta_1, \Theta_2)\sim  \theta$ is no longer necessary. 
To show \eqref{bilinear} in this case,    it suffices by H\"older's inequality to show 
\Be 
\label{small-angle}
\| e^{it\sqrt{-\Delta}} f_j \|_{L^{q}( \mathbb R^d\times  I_\circ)} \lesssim  N^{\frac{d+1}2(\frac1p-\frac{1}q)}   \|f_j\|_{L^p}, \quad j=1,2
\Ee
for $2\le p\le q\le \infty$. The estimate for $(p,q)=(2,2)$ is immediate from Plancherel's theorem. Since $\supp {\widehat f_j}$ is contained in a box of
measure $N^{\frac{d+1}2}$, the Cauchy--Schwarz inequality yields  $\| e^{it\sqrt{-\Delta}} f_j \|_{L^{\infty}( \mathbb R^d\times  I_\circ)} \lesssim N^{\frac{d+1}4} \|\widehat f_j\|_2$.  Thus, the estimate for $(p,q)=(2,\infty)$ follows again by  Plancherel's theorem.  The case $(p,q)=(\infty, \infty)$ is clear from the previous computation at the beginning of the proof. Finally, interpolation between these estimates yields all the  desired estimates for $2\le p\le q\le \infty$.
\end{proof}

\newcommand{\diam}{\operatorname{diam}}

{By} finite angular decomposition and rotation,  we may assume  
\Be 
\label{supp-f} \supp \widehat f\subset  \mathbb A_N^0:= \{ \xi:  \xi/|\xi|\in \Theta_0, \ N/2\le  |\xi|\le 2N \},
\Ee
where $\Theta_0$ is a small spherical cap around  $e_1$. 
We make a Whitney type decomposition of  $\Theta_0\times {\Theta_0}$  away from its diagonal $\{ (\theta, \theta): \theta\in \Theta_0
\}$.  

\subsubsection*{Whitney type decomposition.} 
Let  $\nu_\circ=\nu_\circ(N)$ denote 
 the integer  $\nu_\circ$ such that 
\[  N^{-1} <   2^{-2\nu_\circ} \le  4N^{-1} .  
\] Following a typical dyadic decomposition process, for each positive integer $0\le \nu \le  \nu_0$, we partition
$\Theta_0$ into spherical caps $\Theta_k^\nu$  such that  
\[c_d  2^{-\nu}\le 
  \diam (\Theta_k^\nu)  \le C_d2^{-\nu}\] for some constants $c_d$, $C_d>0$
and 
$\Theta_k^\nu \subset \Theta_{k'}^{\nu'}$ for some $k'$ whenever $\nu\ge \nu'$. 
Thus, we have
\Be
\label{thetak} \Theta_0=\bigcup_{k}\Theta_k^\nu\Ee 
for each $\nu$. Consequently,  we may write
\Be
\label{bi-decomp} \Theta_0\times \Theta_0=\bigcup_{0\le \nu\le \nu_\circ }\,\, \bigcup_{k\sim_\nu k'} \Theta_k^\nu\times \Theta_{k'}^{\nu},\Ee
where  $k\sim_\nu k'$ means $\dist (\Theta_k^{\nu}, \Theta_{k'}^{\nu})\sim 2^{-\nu}$ if  $\nu< \nu_\circ$ and  $\dist (\Theta_k^{\nu}, \Theta_{k'}^{\nu})\lesssim 2^{-\nu}$ if  $\nu= \nu_\circ $ 
(e.g., see \cite[p.\,971]{TVV}). When  $k\sim_{\nu_\circ} k'$, the sets $\Theta_k^{\nu_\circ}$ and $\Theta_{k'}^{\nu_\circ}$ are not necessarily separated from each other by 
distance $\sim 2^{-\nu_\circ}$ since  the decomposition 
process terminates  at $\nu=\nu_\circ$.  

We now consider $L^p\times L^p\to L^{q/2}$ estimate for the bilinear operator $(f, g)\mapsto e^{it\sqrt{-\Delta}} f  e^{it\sqrt{-\Delta}} g$,  whose boundedness  clearly implies that of the linear operator $f\mapsto e^{it\sqrt{-\Delta}} f$ if we set $g=f$. 
We assume that 
\Be
\label{supp-fg}   \supp \widehat f\subset \mathbb A_N^0, \quad  \supp \widehat g \subset  \mathbb A_N^0.\
\Ee

For each $\nu$, let  $\{\tilde\chi_{\Theta_k^\nu} \}_{k}$ be a partition of unity subordinated to $\{\Theta_k^\nu\}_{k}$ with 
$\partial^\alpha \tilde\chi_{\Theta_k^\nu}=O(2^{\nu|\alpha|})$.  Let us define $f_k^\nu$ and $g_{k'}^\nu$ by
\[  \mathcal F (f_k^\nu)( \xi)=    \widehat f(\xi) \tilde\chi_{\Theta_k^\nu} \Big(\frac{\xi}{|\xi|}\Big), \quad    \mathcal F (g_{k'}^\nu)( \xi)=    \widehat g(\xi) \tilde\chi_{\Theta_{k'}^\nu} \Big(\frac{\xi}{|\xi|}\Big)     .\] 
By \eqref{thetak} and  \eqref{supp-fg}, it is clear that $f=\sum_k f_k^\nu$ and $g=\sum_k g_k^\nu$ for each $\nu$. Therefore, we have
 \Be
\label{bi-decomp1} e^{it\sqrt{-\Delta}} f e^{it\sqrt{-\Delta}} g= \sum_{{0\le \nu \leq \nu_\circ}} \sum_{k\sim_\nu k'} e^{it\sqrt{-\Delta}}  f_k^\nu(x)\, e^{it\sqrt{-\Delta}} g_{k'}^\nu(x).
\Ee

Now, for each $\nu$,  define a bilinear operator by 
\[  B^\nu (f,g) (x,t)=   \sum_{k\sim_\nu k'} e^{it\sqrt{-\Delta}}  f_k^\nu(x)\, e^{it\sqrt{-\Delta}} g_{k'}^\nu(x).   \] 
From \eqref{bi-decomp1}, it  follows  that 
\Be \label{bi-decom}   e^{it\sqrt{-\Delta}} f  e^{it\sqrt{-\Delta}} g =\sum_{{0\le \nu \leq \nu_\circ}}    B^\nu (f,g).\Ee

Thanks to orthogonality of the family  $\{e^{it\sqrt{-\Delta}}  f_k^\nu e^{it\sqrt{-\Delta}} g_{k'}^\nu\}_{{k\sim_\nu k'}}$,  
the task of estimating  $L^p \times L^p\to L^{q/2}(\mathbb{R}^d \times \mathcal E)$ operator norm of $B^\nu(f, g)$ reduces to 
 that of  $(f, g)\mapsto e^{it\sqrt{-\Delta}} f_k^\nu e^{it\sqrt{-\Delta}} g_{k'}^\nu$ with, $k\sim_\nu k'$.

\begin{lem} 
\label{bisum}  Let  $1/2 \le 1/p+1/q$ and $2\le p\le q$.  Let  $0\le \nu\le \nu_\circ$. Suppose that 
\Be
\label{assum}
 \|e^{it\sqrt{-\Delta}}  f_k^\nu\, e^{it\sqrt{-\Delta}} g_{k'}^\nu \|_{L^{q/2}( \mathbb R^d\times \mathcal E)}\le \mathcal B \|f_k^\nu\|_p\|g_{k'}^\nu\|_p
 \Ee
 holds for a constant $\mathcal B$ whenever $k\sim_\nu k'$.
Then, there is a constant   $C>0$, independent of $\nu$, such that
\[\|B^\nu (f,g) \|_{L^{q/2}( \mathbb R^d\times \mathcal E)}\le  C\mathcal B \|f\|_p\|g\|_p.\]
\end{lem}

To prove the above lemma, we need the following result, which is basically the same as \cite[Lemma 7.1]{TV2} but in a slightly different form. For clarity we provide its proof below.

\begin{lem}\label{TVV-lemma} Let $r\in [1,\infty]$.  Let $\{R_k\}_{k}$ be a collection of rectangles such that $\{2R_k\}_k$
  overlap at most $M$ times (here, $2R_k$  denotes the dilation of $R_k$  by a factor $2$ with respect to its center). Suppose $\supp \widehat {F_k} (\cdot, t)\subset R_k$ for all $t\in \mathcal E$.  Then,  we have 
 \[  \| \sum_k F_k  \|_{L^r( \mathbb R^d\times \mathcal E)} \le M^{\min(\frac1r, \frac1{r'})}   \Big( \sum_k  \| F_k  \|_{L^r( \mathbb R^d\times \mathcal E)}^{\min(r, r')}\Big)^{1/\min(r, r')}.  \]
\end{lem}

\begin{proof} 
Let $ m_k$ be a smooth  Fourier multiplier adapted to $R_k$ in a natural manner  such  that  $\supp m_k\subset 2R_k$ and  $m_k=1$ on $R_k$.
By considering  the vector valued operator $ \{f_k\}_k  \mapsto   \{m_k(D) f_k\}_k$, 
it is sufficient to show 
 \Be 
 \label{otho0} \| \sum_k  m_k(D)  f_k  \|_{L^r( \mathbb R^d\times \mathcal E)}\le  M^{\min(\frac1r, \frac1{r'})}   \Big( \sum_k  \| f_k  \|_{L^r( \mathbb R^d\times \mathcal E)}^{\min(r, r')}\Big)^{1/\min(r, r')}.   \Ee
For a fixed $t\in \mathcal E$,  by Plancherel's theorem it follows that 
\[ \| \sum_k  m_k(D)  f_k(\cdot, t)  \|_{L^2( \mathbb R^d)}^2 \le  M \Big( \sum_k  \| f_k (\cdot, t)  \|_{L^2( \mathbb R^d)}^{2}\Big).\]
Summation over $t\in \mathcal E$ gives \eqref{otho0} with $r=2$.  On the other hand, \eqref{otho0} for $r=1$ and $r=\infty$ is clear by the triangle inequality. 
 Interpolation (complex interpolation) between the consequent estimates gives  \eqref{otho0}  for all $1\le r\le \infty.$ 
\end{proof}

\begin{proof}[Proof of Lemma \ref{bisum}] 
We begin by noting that, for each $t\in \mathbb R$,   
\[\supp \mathcal F_x(\{ e^{it\sqrt{-\Delta}}  f_k^\nu e^{it\sqrt{-\Delta}} g_{k'}^\nu\}_{k\sim_\nu k'})\]
are contained in the rectangles $R_k^\nu$ of dimensions about \[ \textstyle N\times \overbrace{ N2^{-\nu}\times\dots \times N2^{-\nu}}^{(d-1)-\text{times}},\] whose longest sides are  parallel to the centers of $\Theta_k^\nu$. Since the centers of $\Theta_k^\nu$ are separated by  a distance $\sim 2^{-\nu}$, it is clear that  the  rectangles $R_k^\nu$ are overlapping at most a fixed constant number of times.

Thus, by Lemma \ref{TVV-lemma} we have 
\[
\| B^\nu (f,g)\|_{L^{\frac q2}( \mathbb R^d\times  \mathcal E)} 
\lesssim  \Big( \sum_{k\sim_\nu k'} \|e^{it\sqrt{-\Delta}}  f_k^\nu e^{it\sqrt{-\Delta}} g_{k'}^\nu  \|_{L^{\frac q2}( \mathbb R^d\times  \mathcal E)}^{q_\star} \Big)^\frac{1}{q_\star} ,
\]
where $q_\star=\min(q/2, (q/2)')$.  Hence,  using the assumption \eqref{assum} and  Cauchy--Schwarz inequality, we have
\begin{align*}
\| B^\nu (f,g)\|_{L^{\frac q2}( \mathbb R^d\times \mathcal E)} \lesssim   \mathcal B \Big( \sum_{k}  \|f_k^\nu\|_p^{2q_\star}\Big)^\frac{1}{2q_\star}
 \Big( \sum_{k'}   \|g_{k'}^\nu\|_{p}^{2q_\star}\Big)^\frac{1}{2q_\star}. 
\end{align*}

Since $1/2 \le 1/p+1/q$ and $2\le p\le q$,  it follows that  $2\le p\le 2q_\star$.  Thus, to complete the proof,  we only have to show the inequalities 
\[ \Big(\sum_{k}  \|f_k^\nu\|_{L^p}^{p}\Big)^\frac1p \lesssim \|f\|_{L^p}, \quad \Big(\sum_{k'}   \|g_{k'}^\nu\|_{L^p}^{p}\Big)^\frac1p\lesssim \|g\|_{L^p}\]
for $2\le p\le \infty$. This can be done by  using the similar orthogonality argument and interpolation as in the proof of Lemma \ref{TVV-lemma}. Indeed, we only verify 
the inequalities for $p=2$ and $p=\infty$.  Note that  $f=\sum_k f_k^\nu$ and $g=\sum_k g_k^\nu$ for each $\nu$, so the case $p=2$ follows from  the Plancherel's 
theorem since the Fourier supports of $\{ f_k^\nu\}_k$ (also $\{ g_k^\nu\}_k$)  are contained in  boundedly overlapping rectangles.   
The case $p=\infty$, i.e., $\|f_k^\nu\|_\infty \lesssim \|f\|_\infty$ and $\|g_k^\nu\|_\infty \lesssim \|g\|_\infty$,   is clear from a kernel estimate. Indeed, note that $ f_k^\nu=K_k^\nu\ast f $ and $\|K_k^\nu\|_1\le C$ for a constant.  
\end{proof}

\begin{prop}
\label{bi-1}
Let $ \frac{d-1}{d+1}< \alpha \le 1$ and  $1/\sqrt N \le \theta \le 1/10$.   
Let $\mathcal E\subset [1,2]$ be $N^{-1}$-separated and $\mathcal A_\alpha:=\mathcal A_\alpha(\cE; N^{-1})$ be given by \eqref{agamma}.   
Suppose that $\Theta_1, \Theta_2\subset \mathbb S^{d-1}$ are spherical caps of diameter $\theta$ such that $\dist(\Theta_1, \Theta_2)\sim  \theta$ if 
 $\theta> 2 N^{-1/2}$, and  $\dist(\Theta_1, \Theta_2)\le \theta$  if $\theta \le 2 N^{-1/2}$.    
 Suppose  that   $\supp \widehat f_1$ and $\supp \widehat f_2$ satisfy  \eqref{supp-con}.   Then, we have 
\Be 
\label{bilineardiscrete} 
\| e^{it\sqrt{-\Delta}} f_1 e^{it\sqrt{-\Delta}} f_2 \|_{L^{q /2}( \mathbb R^d\times \mathcal E)} \le C  B^\alpha_{p,q}(N,\theta) \|f_1\|_{p} \|f_2\|_{p}
\Ee
for $q\ge q_\circ$  and $2/p=q_\circ/q$, where 
\Be 
\label{BN}  
B^\alpha_{p,q}(N,\theta)
= \mathcal A_{\alpha}^\frac2{q}  N^{2 (\frac{d-1}2+\frac{1}{p}- \frac{d}q)}
 \theta^{ 2((d-1)(1-\frac1p)- \frac{d-1+2\alpha}{q})}. 
 \Ee 
\end{prop}

\begin{proof}
The proof of  \eqref{bilineardiscrete} is based on interpolation between the following:
\begin{align} 
\label{bilineardiscrete0} 
\| e^{it\sqrt{-\Delta}} f_1 e^{it\sqrt{-\Delta}} f_2\|_{L^{\infty}( \mathbb R^d\times \mathcal E)} &\lesssim  
B^\alpha_{\infty,\infty}(N,\theta) \|f_1\|_{\infty} \|f_2\|_{\infty}, 
\\
\label{bi-discrete2} 
\| e^{it\sqrt{-\Delta}} f_1 e^{it\sqrt{-\Delta}} f_2 \|_{L^{q_\circ /2}( \mathbb R^d\times \mathcal E)} &\lesssim  
B^\alpha_{2,q_\circ}(N,\theta)
 \|f_1\|_{2} \|f_2\|_{2}.
\end{align}
The first follows from \eqref{bilinear} with $p=q=\infty$.  Thus it remains to show \eqref{bi-discrete2}. To this end, we make use of the locally constant property related to frequency localization.    

As before, by rotational symmetry, we may assume that  \eqref{supp0} holds. Let us set  
\Be 
\label{ut}  U_t h= \int e^{i(x\cdot\xi+ t\phi(\xi) )} \widehat h(\xi) \,d\xi,
\Ee
where $\phi(\xi)=|\xi|-\xi_1$. 
Making change of variables  $x_1\to x_1-t$, we note that
\Be
\label{ut-lq}
\| e^{it\sqrt{-\Delta}} f_1 e^{it\sqrt{-\Delta}} f_2 \|_{L^{q}( \mathbb R^d)} = \| U_t f_1  U_t f_2 \|_{L^{q}( \mathbb R^d)}.
\Ee
 Thus,  the desired estimate \eqref{bi-discrete2} is equivalent to
\Be 
\label{bilineardiscrete2q02} 
\| U_t f_1  U_t f_2 \|_{L^{q_\circ /2}( \mathbb R^d\times \mathcal E)} \le C \mathcal A_{\alpha}^\frac2{q_\circ} N^{2 (\frac{d}2- \frac{d}{q_\circ})}
 \theta^{ 2(\frac{d-1}2- \frac{d-1+2\alpha}{{q_\circ}})} \|f_1\|_{2} \|f_2\|_{2}.
 \Ee
 
Observe that 
\Be 
\label{trans-phi}
 \phi(\xi)=  \xi_1 \tilde \phi( \xi'/\xi_1),  \quad \tilde \phi(\eta):=(1+|\eta|^2)^{1/2}-1. 
 \Ee
  Since  $\tilde \phi(\eta)=\frac12 |\eta|^2 +O(|\eta|^4)$, recalling \eqref{supp0}, we  note that $\phi(\xi)=O(N\theta^2)$ for $\xi\in  \supp\widehat f_1\cup \supp \widehat f_2$. Consequently, 
  if $|t-t'|\le N^{-1} \theta^{-2}$,     we may heuristically treat  $U_t f_j$ and $U_{t'}f_j$ 
as if  
\[
 |U_t f_j(x)| \sim |U_{t'}f_j(x)|. 
\] 
 This can be made rigorous making use of the Fourier series expansion as before (cf., Lemma \ref{lcp}).   

 Let $\{I\}$ be a  collection of  almost disjoint intervals of length $\sim N^{-1} \theta^{-2}$ that cover $\mathcal E$ and satisfy $\mathcal E\cap I\neq \emptyset$ for every $I$. 
 For each interval $I$, let $t_I$ denote the center of $I$. 
   Then,   using \eqref{supp0} and   the same argument as in the proof of Lemma \ref{lcp}, we obtain  
\Be
\label{lc-theta}
 |U_t f_j (x)|  \le     \mathfrak S(U_{t_I} f_j )(x)\le  \mathfrak S^2(U_{t} f_j )(x) , \quad   j=1,2
\Ee
for  all $t\in I$. Consequently,  it follows that 
\begin{align*}
   \sum_{t\in \mathcal E} |U_t f_1 (x) U_t f_2 (x)|^q &\lesssim  \sum_{I} \#(\mathcal E\cap I) \big( \mathfrak S(U_{t_I} f_1 )(x) \mathfrak S( U_{t_I} f_2)(x)\big)^q 
   \\ 
   & \lesssim   \sigma  \sum_{I}  \big( \mathfrak S(U_{t_I} f_1 )(x) \mathfrak S( U_{t_I} f_2)(x)\big)^q
   \end{align*}
for any $1\le q<\infty$, where
\Be
\label{equ-sig}
 \sigma=  \sup_{I: |I|\sim  N^{-1}\theta^{-2}}  \#(\mathcal E \cap I). 
 \Ee  
 Since the length of the intervals $I$ is  $\sim  N^{-1} \theta^{-2}$,  the second inequality in \eqref{lc-theta} leads to 
\[  \sum_{I}  \big( \mathfrak S(U_{t_I} f_1 )(x) \mathfrak S( U_{t_I} f_2)(x)\big)^q  \lesssim  N \theta^{2}   \int_{[1,2]}  \big( \mathfrak S^2(U_{t} f_1 )(x) \mathfrak S^2( U_{t} f_2)(x)\big)^q  dt
   \] 
   for any $1\le q<\infty$. 
Combining the inequalities together and taking integration in $x$,  we particularly have
 \Be
  \label{hoho} \| U_t f_1  U_t f_2 \|_{L^{q_\circ /2}( \mathbb R^d\times \mathcal E)} \lesssim (N \theta^{2}\sigma)^{2/q_\circ} \| \mathfrak S^2(U_{t} f_1 )  \mathfrak S^2(U_{t} f_2)\|_{L^{q_\circ/2}(\mathbb R^d\times [1,2])}. 
  \Ee

We now recall \eqref{def-sj}.  Thanks to the rapid decay of $c_{\mathbf k}$, it is  clear that   $L^2\times L^2\to L^{q_\circ/2}(\mathbb R^d\times [1,2])$   norm $\mathfrak N$ of  $(f,g)\mapsto 
 \mathfrak S^2(U_{t} f )  \mathfrak S^2(U_{t} g)$  is bounded by  a constant times 
 that of $(f,g)\mapsto U_{t}f  U_{t} g$ (see Remark \ref{rem:s-sum}).  Moreover,  it is  easy to see  
 \[ \| U_{t}f  U_{t} g\|_{L^{q_\circ/2}(\mathbb R^d\times [1,2])}=\| e^{it\sqrt{-\Delta}} f  e^{it\sqrt{-\Delta}} g \|_{L^{q_\circ/2}( \mathbb R^d\times [1,2])}.\] 
 Thus, using 
\eqref{bilinear} with $q=q_\circ$ and $p=2$, we obtain  
 \[ \mathfrak N\lesssim  N^{2 (\frac{d}2-  \frac{d+1}{q_\circ})}
 \theta^{ 2(\frac{d-1}2- \frac{d+1}{q_\circ})}.\]
 
 Consequently, combining this and \eqref{hoho}  yields  
\[\| e^{it\sqrt{-\Delta}} f_1  e^{it\sqrt{-\Delta}} f_2 \|_{L^{q_\circ /2}( \mathbb R^d\times \mathcal E)} \le CN^{2 (\frac{d}2- \frac{d+1}{q_\circ})}
 \theta^{ 2(\frac{d-1}2- \frac{d+1}{q_\circ})}(N \theta^2 \sigma)^{\frac2{q_\circ}}  \|f_1\|_{2} \|f_2\|_{2}.\]
From \eqref{agamma}, we  also have  
\Be \label{sig} \sigma\le \mathcal A_{\alpha} \theta^{-2\alpha}.\Ee
 Therefore, \eqref{bilineardiscrete2q02} follows as desired. 
\end{proof} 

We  are now ready  to prove Proposition \ref{p-ls}.

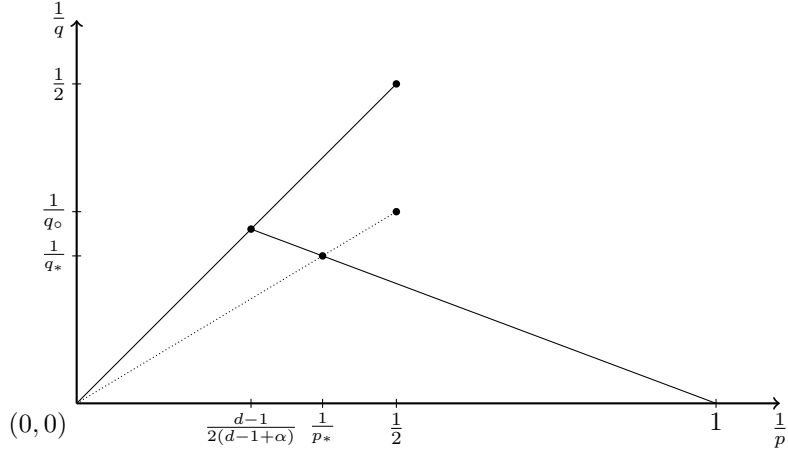
\begin{figure}[t] 
	\centering 
	\begin{tikzpicture}[scale=0.85]
		
		\draw[thick,->]  (0,0) -- (11,0);     
		\draw[thick,->] (0,0) -- (0,6);
		\node[below] at (11,0) {$\frac{1}{p}$};
		\node[left] at (0,6) {$\frac{1}{q}$};
		
		\node[below left] at (0,0) {$(0,0)$};  
		
		\draw (10,-2pt) -- (10,2pt);           
		\node[below] at (10,0) {$1$};          

		\draw (-2pt,5) -- (2pt,5);  
		\node[below] at (5,0) {$\frac{1}{2}$};
		
		\draw (5,-2pt) -- (5,2pt);             
		\node[left] at (0,5) {$\frac{1}{2}$}; 
		
		\draw (-2pt,3) -- (2pt,3);        
		\node[left] at (0,3) {$\frac{1}{q_{\circ}}$}; 
		
        \draw (50/13,-2pt) -- (50/13,2pt);      
        \node[below,font=\footnotesize] at (50/13,0) { $\frac{1}{p_{\ast}}$};
        
        \draw (-2pt,30/13) -- (2pt,30/13);      
        \node[left,font=\footnotesize] at (0,30/13) { $\frac{1}{q_{\ast}}$};

		\node[below,font=\footnotesize] at (30/11,0) { $\frac{d-1}{2(d-1+\alpha)}$};

		\node[circle, fill=black, inner sep=1.0pt] at (5,3) {};
		\node[circle, fill=black, inner sep=1.0pt] at (5,5) {};
		\node[circle, fill=black, inner sep=1.0pt] at (50/13,30/13) {};
		\node[circle, fill=black, inner sep=1.0pt] at (30/11,30/11) {};
		\draw (30/11,-2pt) -- (30/11,2pt);
		
		\draw (0,0) -- (5,5);     
		\draw[densely dotted,black] (0,0) -- (5,3);      
		\draw (30/11,30/11) -- (10,0);
	\end{tikzpicture}
	\caption{Treating the $L^2 \times L^2 \to L^{q_\circ/2}$ estimate as a linear $L^2 \to L^{q_\circ}$ bound, and applying Bourgain’s trick with the $L^\infty$ estimate,  yields a restricted weak-type bound at $(p_\ast, q_\ast)$, as illustrated above for $d = 2$ and $\alpha = 5/6$.} 
	\label{Fig2} 
\end{figure}

 \begin{proof}[Proof of Proposition \ref{p-ls}]  
 Recall $q_\ast= \frac{2(d^2+(\alpha+1) d+\alpha-2)}{d^2-1}$ and set 
 \Be\label{pq}  p_\ast=2q_\ast/ {q_\circ}.\Ee
  Since $ \frac{d-1}{d+1}< \alpha$, $q_\ast> q_\circ= 2(d+3)/(d+1)$.
We also note that  $(p,q)=(p_\ast,q_\ast)$ satisfies \eqref{s-line} (see Figure \ref{Fig2}).  
Thus, to prove Proposition \ref{p-ls}, it is sufficient to show  the restricted weak type estimate 
\Be
\label{res-weak} \| e^{it\sqrt{-\Delta}} f\|_{L^{q_\ast, \infty} ( \mathbb R^d\times \mathcal E)} \le C \mathcal A_{\alpha}^\frac1{q_\ast} 2^{(\frac{d-1}2+\frac1{p_\ast}-\frac d{q_\ast}  )j} \|f\|_{p_\ast,1} 
\Ee 
whenever $\widehat f$  is supported in  the annulus $\mathbb A_{2^j}^0$.
Since  the estimate \eqref{localsmoothing-2d1} holds for $p=1$ and $q=\infty$,  applying  real interpolation yields the estimate \eqref{localsmoothing-2d1}  for $q>q_\ast$.

To show \eqref{res-weak}, we take $N=2^j$ and make  use of  \eqref{bi-decom}.  For $f,g$ satisfying \eqref{supp-fg}, we consider the operator $B^\nu (f,g)$. 
  Proposition \ref{bi-1} gives
 \[ 
 \|e^{it\sqrt{-\Delta}}  f_k^\nu e^{it\sqrt{-\Delta}} g_{k'}^\nu  \|_{L^{\frac q2}( \mathbb R^d\times  \mathcal E)}\lesssim  B^\alpha_{p,q}(N, 2^{-\nu})   \|f\|_{p} \|g\|_{p}
\]
for $q\ge q_\circ$ and $2/p=q_\circ/q$ provided that  $k\sim_\nu k'$. Thus, applying Lemma~\ref{bisum},  we obtain 
\Be
\label{bnu}
\|  B^\nu (f,g)\|_{L^{q /2}( \mathbb R^d\times \mathcal E)} \lesssim \mathcal A_{\alpha}^\frac2{q} N^{2 (\frac{d-1}2+\frac{1}{p}- \frac{d}q)}
 2^{ 2\nu (\frac{d-1+2\alpha}{q}-(d-1)(1-\frac1p))}  \|f\|_{p} \|g\|_{p}
\Ee
if  $q\ge q_\circ$, $2/p=q_\circ/q$, and  $(p,q)$ is contained  in a sufficiently small neighborhood of $(p_\ast, q_\ast)$.   Indeed, since $1/2 < 1/p_\ast+1/q_\ast $ and $2< p_\ast< q_\ast$, 
there is a  neighborhood  $V$ of $(p_\ast, q_\ast)$ such that   $1/2\le 1/p+1/q$ and $2\le p\le q$ hold for  $(p,q)\in V$. 
 
Recall  \eqref{bi-decom}.  We  now apply Lemma \ref{B-Sum-Tri} to  the operators $\{B^\nu (f,g)\}_{0\le \nu \le n_\circ},$ using the estimate \eqref{bnu}.  In fact, we only need to choose two pairs $(p_0, q_0)$ and $(p_1, q_1)$ in $V$ such that $q_j \ge q_\circ$, $2/p_j=q_\circ/q_j$ for $j=0,1$, and 
\[ \frac{d-1+2\alpha}{q_0}-(d-1)\Big(1-\frac1{p_0}\Big)>0> \frac{d-1+2\alpha}{q_1}-(d-1)\Big(1-\frac1{p_1}\Big). \]  Note that $1-\frac1{p_\ast}= \frac{d-1+2\alpha}{q_\ast}$. Therefore,  
by Lemma \ref{B-Sum-Tri}  with $n=2$, we obtain 
 \[  \| e^{it\sqrt{-\Delta}} f  e^{it\sqrt{-\Delta}} g\|_{L^{q_\ast /2, \infty}( \mathbb R^d\times \mathcal E)} \le C \mathcal A_{\alpha}^\frac2{q_\ast} 
 2^{2j (\frac{d-1}2+\frac{1}{p_\ast}- \frac{d}{q_\ast})} \|f\|_{p_\ast,1} \|g\|_{p_\ast,1}, \] 
which yields \eqref{res-weak} if we take $ g=f $.  
\end{proof}

\section{Improvement of the smoothing estimate via a bilinear estimate}
\label{sec:improved}

In the previous section, the estimate \eqref{smooth-pq} is  verified  for $q> q_\ast$  and $ s=s_1(p,q) $  while \eqref{assouad} and \eqref{s-line} are satisfied.  Though the range is wider than those in \cite{BRRS, Wh}, the estimate remains subjected to the restriction \[ \alpha>(d-1)/(d+1).\] 
This limitation arises from  our direct use of the bilinear restriction estimate \eqref{bilinear0}, whose range  of $q$ is inherently constrained by the one-dimensional underlying set in the 
$t$-variable.  However, when we consider the sets of fractal dimension for $t$, the underlying set  becomes sparse as the dimension gets smaller. Therefore, it is natural to expect that one can further improve the range of  bilinear estimate by incorporating  the sparseness of the $t$-set $E$.

That is to say, we  extend the  range of $q$ ($p$ given by \eqref{s-line}) for which the estimate \eqref{smooth-pq}  with $ s=s_1(p,q)$ holds to  $q> \tilde q_\ast$, while also removing  the restriction on $\alpha$, where 
\[ \tilde q_\ast:= \frac{2(d-1+2\alpha)^2-4\alpha^2}{(d-1)(d-1+2\alpha)}.\]
This is achieved  by establishing Proposition \ref{p-ls'} below.  
However, since we  utilize an induction on scales argument of iterative nature, a quantitative control of the bound with $\mathcal A_\alpha(\cE; 1/N)$,  such  as the one in Proposition \ref{p-ls} (see \eqref{agamma}  and \eqref{localsmoothing-2d1}), is no longer available.  Moreover,  a stronger condition on the set 
$\cE$ than \eqref{agamma} is required for our purposes.

Recalling \eqref{agamma}, 
we define
\Be
\label{agamma'}
\tilde{\mathcal A}_{\alpha}(\cE; \delta)=    \sup_{\delta \le  \delta' < 1}   {\mathcal A}_{\alpha}(\cE; \delta'). 
\Ee
Though $\tilde{\mathcal A}_{\alpha}$ and  ${\mathcal A}_{\alpha}$ are similar in form,   the difference between them can be significant in some cases. 
For instance, for $0<\beta<\alpha<1$,  consider the $\delta$-separated set $\cE=\{1+k\delta^{\beta}: 1\leq k\leq \delta^{-\beta} \}$.     It is easy to see that $ \tilde{\mathcal A}_{\alpha}(\cE; \delta)\gtrsim \delta^{\beta(\alpha-1)}$ while $ {\mathcal A}_{\alpha}(\cE; \delta)=1$.

The proof of Theorem \ref{improved} essentially reduces to the following proposition.   

\begin{prop}\label{p-ls'} Let $0< \alpha \le 1$,  $N\ge 1$, and   $ \cE $ be a $1/N$-separated subset of $[1,2]$.  
Suppose that $\widehat f$ is supported in $\mathbb A_N$ and  $\tilde{\mathcal A}_{\alpha}(\cE; 1/N)\le \tilde{\mathcal A}$ for a constant  $\tilde{\mathcal A}$.
Then, for $q> \tilde q_\ast$ and $p$ given by \eqref{s-line},    there is a constant $C=C(\tilde{\mathcal A})$ such that 
\Be 
\label{localsmoothing-2d} 
\| e^{it\sqrt{-\Delta}} f\|_{L^q( \mathbb R^d\times  \cE)} \le C N^{\frac{d-1}2+\frac1p-\frac dq} \|f\|_p. 
\Ee    
\end{prop}

Indeed, if $E$ is a set of bounded $\alpha$-Assouad characteristic, then the $2^{-j}$-separated set 
  $\mathcal E=E(2^{-j})$ satisfies $\tilde{\mathcal A}_{\alpha}(\mathcal E ; 2^{-j})\le \mathcal A^\alpha(E)$ for  all $j$. Therefore, Theorem \ref{improved} follows 
  from Proposition \ref{p-ls'}.

Let $\Gamma\subset \mathbb R$.  To prove Proposition \ref{p-ls'} we first  obtain   a bilinear restriction estimate over $\mathbb R^d\times  \Gamma$, where $\Gamma$ is a sparse union of essentially  disjoint unit intervals. 
To make it precise, let   $\mathcal I$ be a collection of disjoint unit intervals $\bf I$ such that 
\Be 
\label{intervals}
\#\{ \mathbf I\in \mathcal I:  \mathbf I  \cap (t, t+r)\neq \emptyset \}  \le  Cr^\alpha 
\Ee
for any $t\in \mathbb R$ and $r\ge 1$.  We then set
\Be  \label{TT} \Gamma =\bigcup_{\mathbf I\in \mathcal I} \, \mathbf I.   \Ee
Recall  the adjoint restriction operator  $\mathcal R^\ast$  defined by \eqref{rest-o}. 

\begin{thm}\label{bilinear-stronger}  Let $\Theta, \Theta'\subset \mathbb S^{d-1}$ with $\dist(\Theta, \Theta')\sim 1$.  Suppose  that  $\mathcal I$ satisfies \eqref{intervals} and $\Gamma$ is given by \eqref{TT}. Then, for $q>\tilde q_\circ:=\frac{2(d-1+4\alpha)}{d-1+2\alpha}$,  we have the estimate   
\Be 
\label{bi-r} \|  \mathcal R^\ast f\, \mathcal R^\ast g \|_{L^{q /2}( \mathbb R^d\times \Gamma)} \le C \|f\|_{2} \|g\|_{2},\Ee
whenever 
\Be\label{supp-f,g}  \supp f \subset \{ \xi \in \mathbb{A}_1:  \xi/|\xi|\in \Theta\}, \quad \supp g\subset  \{ \xi \in\mathbb{A}_1:  \xi/|\xi|\in \Theta'\}.\Ee
\end{thm}

Without the  separation  condition $\dist(\Theta, \Theta')\sim 1$, one can easily see, using the Knapp example, that the estimate \eqref{bi-r}  holds only if 
\Be
\label{l2-est}
 q\ge 2(d-1+2\alpha)/(d-1).\Ee  Also, note  that $\tilde q_\circ< 2(d-1+2\alpha)/(d-1)$  for any $0<\alpha\le 1$. The bilinear estimate under the separation 
 condition has a wider range of boundedness compared with the linear estimate.   The squashed cap example  (see \cite[Example 14.7]{FK})  shows that  the estimate holds only if 
 \Be
 \label{bi-l2} q\ge 2(d+1+2\alpha)/(d+1).
 \Ee

 Before proceeding  to prove Theorem \ref{bilinear-stronger}, we verify the conditions
 \eqref{l2-est} and \eqref{bi-l2}. To this end, we use the set constructed in Lemma \ref{lem:existA}.  Let $0<\delta\ll 1$ and $j$ be the integer such that $ 2^{j-1}\le  \delta^{-2} <2^j.$ 
Let  $\mathcal E$ be the $2^{-j}$ separated set in Lemma \ref{lem:existA} and set  $\bar{\mathcal E}= \{2^j(t-1): t\in \mathcal E\}$.  Thus, $\bar{\mathcal E}$ is $1$-separated and 
contained in $[0, 2^j]$. 
Set 
\[ \Gamma=\bigcup_{t\in  \bar{\mathcal E}} (t, t+1) .\]
Since ${\mathcal A}_{\alpha}(\mathcal E;  2^{-j})\lesssim 1$, it is easy to see that  $\{  (t, t+1) : t\in  \bar{\mathcal E}\}$ satisfies \eqref{intervals}. 

Consider the functions
\[ f=g=  \delta^{(1-d)/2} \chi_{\{\xi\in \mathbb A_1: |\frac{\xi}{|\xi|} -e_1|<\delta\}}.\]  Then $|\mathcal R^\ast f(x,t)|=|\mathcal R^\ast  g(x,t)|\sim  \delta^{(d-1)/2}$ if 
$|x_1+t|\le c$, $|x'|\le c\delta^{-1}$, and $|t|\le c\delta^{-2}$ with a small enough constant $c>0$.   Since there are $\# \bar{\mathcal E}\sim \delta^{-2\alpha},$  the estimate 
\eqref{bi-r} without separation of the supports implies  $ \delta^{d-1}   \delta^{-2(d-1+2\alpha)/q} \lesssim  1.$
Therefore, letting $\delta\to 0$ implies  \eqref{l2-est}.   

Writing $\xi=(\xi_1, \xi')$, we consider the functions 
\[  f= \delta^{-(d+1)/2}\chi_{\{  |\xi_1-1 | \le \delta^2,\, |\xi'|\le \delta   \}}, \quad  g= \delta^{-(d+1)/2}\chi_{\{  |\xi_1+1 | \le \delta^2,\, |\xi'|\le \delta   \}}. \]
We observe that $|\mathcal R^\ast  f(x,t)|, |\mathcal R^\ast  g(x,t)|\sim  \delta^{(d+1)/2}$ if $|x_1|, |t|\le  c\delta^{-2}$ and $|x'|\le c\delta^{-1}$ with a small enough constant $c>0$.
Since $\# \bar{\mathcal E}\sim \delta^{-2\alpha},$
 the estimate 
\eqref{bi-r} gives  $\delta^{d+1}   \delta^{-2(d+1+2\alpha)/q} \lesssim  1.$ Hence, we get  \eqref{bi-l2}.

\subsection{Proof of Theorem \ref{bilinear-stronger}}
For the  proof of Theorem \ref{bilinear-stronger}, we adapt the induction on scales argument which was used to prove sharp bilinear restriction estimates  \cite{W, T, L2}.  

\subsubsection*{Decomposition of $f$ and $g$}   For $R\gg 1$ and $0<\delta\ll 1$, let $Q_R=[-R,R]^{d+1}$. 
Let $\{\mathbf
q\} $ be a collection of cubes of sidelength $\sim R^{1-\delta}$ that partition the cube $Q_R$, so that 
\[ Q_R =\bigcup \mathbf q.\]

\begin{lem}
\label{wavepacket}  Let $R \gg 1$ and $0<\delta\ll 1$.  Let $\Theta, \Theta'\subset \mathbb S^{d-1}$ with  $\dist(\Theta, \Theta')\sim 1$.  Suppose 
\eqref{supp-f,g} holds. 
Then,  for each $\mathbf
q$, $f$ and $g$ can be decomposed such that
\Be\label{decompo}
 f=f_{\mathbf q}+ f_{{\mathbf q}}^{\not\sim} \quad \text{and} \quad  g=g_{\mathbf q}+ g_{{\mathbf q}}^{\not\sim} 
 \Ee
with $f_{\mathbf q}$,
$f_{{\mathbf q}}^{\not\sim}$, $g_{\mathbf q}$,
$g_{{\mathbf q}}^{\not\sim}$ supported in $\mathbb A_1+O(R^{-1/2})$ 
and, for any $0<\epsilon\ll \delta$, the following hold:
\begin{align}
\label{l2sumq}
 &\sum_{\mathbf q} \|f_{\mathbf q}\|_2^2\le CR^\epsilon
\|f\|_2^2,\qquad \sum_{\mathbf q} \|g_{\mathbf q}\|_2^2\le
CR^\epsilon \|g\|_2^2, 
\end{align}
and, with a constant  $c$ independent of $\delta$ and $\epsilon$, 
\Be
\label{bil2} 
\begin{aligned}
\| \mathcal R^\ast  (f_{\mathbf q}) \mathcal R^\ast  (g_{{\mathbf q}}^{\not\sim})\|_{L^2(\mathbf q)}
&+
 \|\mathcal R^\ast  (f_{{\mathbf q}}^{\not\sim} ) R^\ast  (g_{\mathbf q})\|_{L^2(\mathbf q)}
 + 
\|\mathcal R^\ast  (f_{\mathbf q}^{\not\sim}) \mathcal R^\ast  (g_{{\mathbf q}}^{\not\sim})\|_{L^2(\mathbf q)}
\\
&\le C R^{\frac{1-d}4+c\delta+\epsilon} \|f\|_2\|g\|_2.
\end{aligned}
\Ee
\end{lem} 

This lemma,  which appeared  in \cite{CHL}, is essentially a slight rephrasing of \cite[Lemma 3.5]{W}. The statement in \cite{W} is formulated by fixing dyadic scales relevant to various quantifications that arise after performing wave packet decomposition at scale $R^{-1/2}$ (corresponding to spatial scale $R^{1/2}$). Since there are only $O((\log R)^C)$ such dyadic scales after discarding negligible contributions, combining them suffices to show the lemma. This can be seen even more clearly  in \cite{L2}, particularly in the proof of \cite[Theorem 1.3]{L2}, where a simpler argument based on a more direct wave packet decomposition is given.

While it is also possible to prove Theorem~\ref{bilinear-stronger} in the same manner as in \cite{W}, Lemma~\ref{wavepacket} allows us to present a significantly more concise proof.

\begin{proof}[Proof of Theorem \ref{bilinear-stronger}]   Recall   \eqref{TT}. 
We first show that
\Be 
\label{local-R} 
 \|  \mathcal R^\ast f\, \mathcal R^\ast g \|_{L^{\tilde q_\circ /2}( (\mathbb R^d\times \Gamma)\cap Q_R)} \le C R^\gamma \|f\|_{2} \|g\|_{2}
 \Ee
for any $\gamma>0$ with a constant $C=C(\gamma)$ whenever  the collection $\mathcal I$ of intervals   satisfies  \eqref{intervals}.     

 Let 
us denote  by $\mathcal L(\gamma)$ if the inequality \eqref{local-R} holds  with some constant $C=C(\gamma)$ whenever  the collection $\mathcal I$ of intervals   satisfies  \eqref{intervals}. To show  that \eqref{local-R} holds  for any $\gamma>0$, it is sufficient to show
that  there is a $c$, independent of $\delta, \epsilon$, such that 
\Be\label{implies} \mathcal L(\gamma) \implies  \mathcal L(\max\{\gamma(1-\delta)+\epsilon,  c(\delta+\epsilon) \}).\Ee 
Since $\|\mathcal R^\ast f  \mathcal R^\ast g\|_\infty \le |\mathbb A_1|\|f\|_2\|g\|_2$, \eqref{local-R} trivially holds with $\gamma=2(d+1)/\tilde q_0$. Taking small enough $\delta>0$ and $\epsilon$ and  iterating this implication finitely many times, one can show 
\eqref{local-R} for any given $\gamma>0$.

\newcommand{\bq}{\mathbf q}

  Since $\bigcup \mathbf q=Q_R$, using \eqref{decompo},  we have 
\begin{align*}
 \|  \mathcal R^\ast f\, \mathcal R^\ast g \|_{L^{\tilde q_\circ /2}( (\mathbb R^d\times \Gamma)\cap Q_R)} 
 &\le  \sum_{\mathbf q} \|  \mathcal R^\ast f\, \mathcal R^\ast g \|_{L^{\tilde q_\circ /2}( (\mathbb R^d\times \Gamma)\cap \mathbf q)} 
 \\[4pt]
 &\le \rm I_0+  \rm I_1+ \rm I_2+ \rm I_3, 
\end{align*}
where 
\begin{align*}
\mathrm I_0&= \sum_{\mathbf q}   \| \mathcal R^\ast  (f_{\mathbf q})  \mathcal R^\ast  (g_{{\mathbf q}})\|_{L^{\tilde q_\circ/2}((\mathbb R^d\times \Gamma) \cap \mathbf q)}, 
\\
\mathrm I_1 &=\sum_{\mathbf q}  \| \mathcal R^\ast  (f_{\mathbf q}) \mathcal R^\ast  (g_{{\mathbf q}}^{\not\sim})\|_{L^{\tilde q_\circ/2}((\mathbb R^d\times \Gamma) \cap \mathbf q)},
\\
\mathrm I_2 &= \sum_{\mathbf q}
 \|\mathcal R^\ast  (f_{{\mathbf q}}^{\not\sim} ) R^\ast  (g_{\mathbf q})\|_{L^{\tilde q_\circ/2}((\mathbb R^d\times \Gamma) \cap \mathbf q)}, 
\\
\mathrm I_3 &= \sum_{\mathbf q}
\|\mathcal R^\ast  (f_{\mathbf q}^{\not\sim}) \mathcal R^\ast  (g_{{\mathbf q}}^{\not\sim})\|_{L^{\tilde q_\circ/2}((\mathbb R^d\times \Gamma) \cap \mathbf q)}.
\end{align*}

For $\rm I_0$, we use the induction assumption \eqref{local-R}, which gives  
\[ \| \mathcal R^\ast  (f_{\mathbf q})  \mathcal R^\ast  (g_{{\mathbf q}})\|_{L^{\tilde q_\circ/2}((\mathbb R^d\times \Gamma) \cap \mathbf q)}  \le C R^{\gamma(1-\delta)} \sum_\bq \|f_\bq\|_{2} \|g_\bq\|_{2}\] since $\bq$ is of side length $R^{1-\delta}$. Thus, by the Cauchy--Schwarz inequality 
and \eqref{l2sumq}  in Lemma \ref{wavepacket}  we  obtain 
\Be 
\label{i0} \mathrm I_0 \le C R^{\gamma(1-\delta)} \sum_\bq \|f_\bq\|_{2} \|g_\bq\|_{2} \le  C R^{\gamma(1-\delta)+\epsilon} \|f\|_2\|g\|_2. \Ee

For $\mathrm I_j$, $j=1,2,3$,  recalling \eqref{intervals} we first notice from Plancherel's theorem that 
\[ \| \mathcal R^\ast  (h )\|_{L^2((\mathbb R^d\times \Gamma)\cap Q_r)} \le C  r^{\alpha/2} \|h\|_2 \] 
for any cube $Q_r$ of side length $r\ge 1$. Applying  this along with \eqref{l2sumq}, we have 
\[  \|\mathcal R^\ast  (f_{{\mathbf q}})\|_{L^2((\mathbb R^d\times \Gamma)\cap \bq)}, \, \|\mathcal R^\ast  (f_{{\mathbf q}}^{\not\sim} )\|_{L^2((\mathbb R^d\times \Gamma)\cap \bq)} \leq CR^{\alpha(1-\delta)/2}R^{\epsilon/2}\|f\|_2.\] 
Similarly, we also have 
\[ \|\mathcal R^\ast  (g_{{\mathbf q}})\|_{L^2((\mathbb R^d\times \Gamma)\cap \bq)},\,  \|\mathcal R^\ast  (g_{{\mathbf q}}^{\not\sim} )\|_{L^2((\mathbb R^d\times \Gamma)\cap \bq)} \le CR^{\alpha(1-\delta)/2}R^{\epsilon/2} \|g\|_2.\] In particular, 
by the Cauchy--Schwarz inequality  we have
\[ \|\mathcal R^\ast  (f_{{\mathbf q}}) \mathcal R^\ast  (g_{\mathbf q}^{\not\sim} )\|_{L^{1}((\mathbb R^d\times \Gamma)\cap \bq)}  \le CR^{\alpha+\epsilon}  \| f\|_2\|g\|_2 .\]

We now combine this  and 
\[ \| \mathcal R^\ast  (f_{\mathbf q}) \mathcal R^\ast  (g_{{\mathbf q}}^{\not\sim})\|_{L^2((\mathbb R^d\times \Gamma)\cap \bq))}\le C R^{\frac{1-d}4+c\delta+\epsilon} \|f\|_2\|g\|_2,\] 
which follows from  \eqref{bil2}. Note that 
\[  \frac 2{\tilde q_\circ}=  \frac12  \frac{4\alpha}{d-1+4\alpha}+ 1 \frac{d-1}{d-1+4\alpha} .\]
Thus, by H\"older's inequality we have
\[
\|\mathcal R^\ast  (f_{{\mathbf q}}) \mathcal R^\ast  (g_{\mathbf q}^{\not\sim} )\|_{L^{\tilde q_\circ/2}((\mathbb R^d\times \Gamma) \cap \mathbf q)}  \le   C  R^{c(\delta+ \epsilon) } \| f\|_2\|g\|_2   \]
for some positive constants $c$. By the same argument, we obtain  the same estimates for $\mathcal R^\ast  (f_{{\mathbf q}}^{\not\sim} ) \mathcal R^\ast  (g_{\mathbf q})$ and $\mathcal R^\ast  (f_{{\mathbf q}}^{\not\sim} ) \mathcal R^\ast  (g_{\mathbf q}^{\not\sim})$.  Therefore, it follows that 
\[ \begin{aligned}
\| \mathcal R^\ast  (f_{\mathbf q})    \mathcal R^\ast  (g_{{\mathbf q}}^{\not\sim})& \|_{L^{\tilde q_\circ/2}((\mathbb R^d\times \Gamma) \cap \mathbf q)}
+
 \|\mathcal R^\ast  (f_{{\mathbf q}}^{\not\sim} ) \mathcal R^\ast  (g_{\mathbf q})\|_{L^{\tilde q_\circ/2}((\mathbb R^d\times \Gamma) \cap \mathbf q)}
 \\[3pt]
&+ 
\|\mathcal R^\ast  (f_{\mathbf q}^{\not\sim}) \mathcal R^\ast  (g_{{\mathbf q}}^{\not\sim})\|_{L^{\tilde q_\circ/2}((\mathbb R^d\times \Gamma) \cap \mathbf q)}
 \le C  R^{c(\delta+ \epsilon) } \|f\|_2\|g\|_2.
\end{aligned}
\]

Since  there
are as many as $O(R^{(d+1)\delta})$ $\bq$,  taking sum over $\bq$, we obtain the estimate  
\[ \mathrm I_1+  \mathrm I_2+  \mathrm I_3 \le C R^{c(\delta+ \epsilon) } \|f\|_2\|g\|_2 \]
for some $c>0$. Therefore, combining this and \eqref{i0}, we have 
\[ \|  \mathcal R^\ast f\, \mathcal R^\ast g \|_{L^{\tilde q_\circ /2}( (\mathbb R^d\times \Gamma)\cap Q_R)} \le C (R^{(1-\delta)\gamma+\epsilon}+ R^{c(\delta+ \epsilon) })  \|f\|_{2} \|g\|_{2}.\]
This establishes \eqref{implies}.

Since \eqref{local-R} holds for any $\gamma>0$, making use of the $\epsilon$-removal lemma, \cite[Lemma A.3]{BG}\footnote{Although the lemma \cite[Lemma A.3]{BG} is stated with a specific exponent $q=4$ in bilinear setting, it is easy to check that the lemma actually works for any exponent $q>2$.}, we get  
 the desired estimates for $q>\tilde q_\circ$.   For the reader's convenience we add some details.   Assuming $\|f\|_2=\|g\|_2=1$,  
 we need  to show  \[\| \chi_{\mathbb R^d\times \Gamma} \mathcal R^\ast f\, \mathcal R^\ast g \|_{L^{q/2}_{x,t}}\le C.\]
for some $C>0$ if $q>\tilde q_\circ$.  Since $\|  \mathcal R^\ast f\, \mathcal R^\ast g \|_\infty \le |\mathbb A_1| $, this follows if we show 
\Be 
\label{fk}  |F_k| \le C_\epsilon  2^{\epsilon k} 2^{\tilde q_0 k/2 } 
\Ee
for any  $\epsilon>0$ where 
\[  F_k=\{ (x,t)\in \mathbb R^d\times \Gamma:   2^{-k}\le  | (\mathcal R^\ast f  \mathcal R^\ast g)  (x,t)| <2^{-k+1}\}.\]

The inequality \eqref{fk}  can be deduced  from the estimate 
 \Be
 \label{re-} \| (\chi_{\mathbb R^d\times \Gamma} \mathcal R^\ast h_1  \mathcal R^\ast h_2  ) \chi_{E}  \|_{L^{\tilde q_\circ/2}_{x,t}} \le C_\epsilon |E|^\epsilon\Ee
  for any $\epsilon>0$ whenever   $E$ is a union of   $c$-cubes\footnote{Here, $c$-cubes  are  cubes of side length $c$ and $c$ is a fixed constant chosen sufficiently small.}  and $\|h_1\|_2=\|h_2\|_2=1$. This can be obtained by following the argument given  in   \cite{BG} (specifically, the proof of  \cite[Lemma A.3]{BG}), where a covering lemma due to Tao (\cite[Lemma 3.3]{T1}, \cite[Lemma 4.3]{T2}) plays a key role. By Tao's lemma one can cover $E$ with $C|E|^\delta/\delta$ many sparse collection of balls of radius $O(|E|^{C^{1/\delta}})$  for any $\delta>0$ (see \cite[Definition 3.1]{T1} for the precise definition of a sparse collection). By this, the set $E$ can be further assumed to be a union of balls in a sparse collection.  
 Then, the inequality \eqref{re-} follows from a routine adaptation of the argument in \cite{BG}.  The only difference  is that there is the  additional  $ \chi_{\mathbb R^d\times \Gamma}$  factor. However, the presence of this additional factor does not cause any problem since all the relevant quantities are controlled only  in terms of $|E|$. 

We conclude the proof by   showing  how one can obtain   \eqref{fk} from \eqref{re-}. This was not made entirely clear in \cite{BG}. Instead of \eqref{fk},  it is clearly sufficient to show 
\Be \label{Fk}  |F_k\cap B(0,M)| \le C_\epsilon  2^{\epsilon k} 2^{\tilde q_0 k/2 }\Ee
for any $M\gg 1$. Let $a\in C_c^\infty (B(0,2^2))$ such that  $a=1$ on $\mathbb A_1$.   As before, expanding $e^{i((x,t)-(x',t'))(\xi, |\xi|)} a(\xi)$  into the Fourier series in $\xi$, we 
note that  
\Be
\label{domi}
\begin{aligned}
    |\mathcal R^\ast f (x',t')|  & \le     \mathfrak Gf(x,t):=   \mathfrak S\mathcal R^\ast f (\cdot ,t )(x) ,
    \\
    |\mathcal R^\ast g (x',t')|   &\le     \mathfrak G g(x,t):=   \mathfrak S\mathcal R^\ast g (\cdot ,t )(x) 
\end{aligned}
\Ee
provided that $|(x,t)-(x',t')|\le  c\sqrt d$ for a small fixed constant $c>0$.  Also, since $\mathcal R^\ast h (x+\tau_{\mathbf k}, t )= \mathcal R^\ast (h e^{i(\cdot)\cdot\tau_{\mathbf k}})$,    making use of rapid decay of the coefficients $c_{\mathbf k}$,  it is easy to see that \eqref{re-} implies 
\Be
 \label{re-'} \| (\chi_{\mathbb R^d\times \Gamma}  \mathfrak G h_1 \mathfrak G h_2 ) \chi_{E}  \|_{L^{\tilde q_\circ/2}_{x,t}} \le C_\epsilon |E|^\epsilon\Ee
for any $\epsilon>0$ whenever 
 $E$ is a union of $c$-cubes and $\|h_1\|_2=\|h_2\|_2=1$.

 Let $Q(x,t)$  denote the open cube of  side length $c$ that is centered at $(x,t)$. Let us set 
 \[ E=\bigcup_{(x,t)\in F_k \cap B(0,M)} Q(x,t).\]
  Since $F_k \cap B(0,M)\subset E$,   it is sufficient for \eqref{Fk} to  show 
 \[   | E|   \le C_\epsilon  2^{\epsilon k} 2^{\tilde q_0 k/2}.\] 
  From  \eqref{domi}  
 we observe that $ \mathfrak Gf(x,t) \mathfrak Gg(x,t) \ge 2^{-k}$ for $(x,t)\in E$. Applying \eqref{re-'} with $h_1=f$, $h_2=g$, and (a union of $c$-cubes) $E$, we have 
 \[ 2^{-k} |E|^{2/\tilde q_0} \le C_\epsilon \| (\chi_{\mathbb R^d\times \Gamma}  \mathfrak G f \mathfrak G g) \chi_{E}  \|_{L^{\tilde q_\circ/2}_{x,t}} \le  C_\epsilon |E|^\epsilon  \]
 for any $\epsilon>0.$  This gives the desired inequality since $|E|$ is finite.  
\end{proof}

\begin{rem}
\label{ee}
Of course, the similar bilinear estimate as in Theorem \ref{bilinear-stronger} holds for  $\mathcal T$ defined by \eqref{extension} under the condition \eqref{sep-scaled} with $\theta\sim 1$. 
This is clear since   Lemma \ref{wavepacket}  remains valid for $\mathcal T$ (for example, see \cite[Proof of Theorem1.3]{L2}).  
Once we have the lemma, the rest of proof works without modification. 
\end{rem}

Recalling \eqref{trans-phi},   for $\theta\ge 0$,  we set  $\tilde \phi_0(\eta)=\frac12 |\eta|^2$ and 
 \Be  \label{phi-theta} 
 \tilde  \phi_\theta(\eta)=  \theta^{-2}  \tilde  \phi(\theta \eta), \quad \theta>0.
 \Ee 
 Then, we consider the operators 
  \Be
 \label{R*-theta}  \mathcal T_\theta h(x,t)= \int e^{i(x\cdot\xi+ t\xi_1\tilde \phi_\theta(\xi'/\xi_1) )} h(\xi) \,d\xi. \Ee
  In order to prove Proposition \ref{p-ls'},  we use rescaling, which gives  rise to a family of extension operators  $\mathcal T_\theta$, $\theta>0$.  
 This occurs because,  unlike in the proof of Proposition \ref{thm:bi}, the transformation $(x,t)\to (2^{-1/2}(x_1-t), x', 2^{-1/2}(x_1+t))$ (see below \eqref{extension}) is no longer permissible when we are dealing with estimates over a discrete time set.

 For this reason,  we need a uniform estimate for   the operators $\mathcal T_\theta$.  

\begin{cor} \label{bilinear-theta} Let   $q>\tilde q_\circ=2(d-1+4\alpha)/(d-1+2\alpha)$.  Suppose  that  $\mathcal I$ satisfies \eqref{intervals} and $\Gamma$ is given by \eqref{TT}.    If $0\le \theta\le \theta_0$ for a positive  constant  $\theta_0$ small enough,  there is a constant $C$, independent of $\theta$, such that 
\Be
\label{theta-unif}
 \| \mathcal T_\theta h_1\mathcal T_\theta  h_2 \|_{L^{q /2}( \mathbb R^d\times \Gamma)} \le C \|h_1\|_{2} \|h_2\|_{2}
\Ee
whenever   $  \supp h_1 \subset \{ ({\rho,\eta}) : \rho\sim 1,   |\eta| {\sim 1}    \}$ and $\supp h_2\subset \{ ({\rho,\eta}) : \rho\sim 1,  |\eta|  \le c   \}$ for a small positive constant $c$.    
\end{cor}

In fact, this can also be shown without assuming  the smallness of $\theta_0$ by performing additional finite decomposition of $\supp h_1$, $\supp h_2$ into sufficiently small angular sectors.   

Corollary \ref{bilinear-theta} can be verified by showing  that Lemma \ref{wavepacket} and the $\epsilon$-removal argument  remain valid uniformly in $\theta$. 
Once this is achieved, the rest of argument is identical to the proof  of Theorem \ref{bilinear-stronger} . Since  the conic surface $(\xi_1, \xi', \xi_1\tilde \phi_\theta(\xi'/\xi_1))$ is a small perturbation of $(\xi_1, \xi', \xi_1\tilde \phi_0(\xi'/\xi_1))$, this can be done by going through the proofs in \cite{W, L2, BG}. Instead of reproducing all  the details, we provide a brief explanation below. 

 Note that $\mathcal T_0$ coincides with $\mathcal T$, which is given by \eqref{extension}. 
As mentioned in Remark \ref{ee}, the estimate \eqref{bi-r} holds for $\mathcal T_0$.  Since $\tilde \phi(\eta)=\frac12 \eta^2 +O(\eta^4)$, we have
\[\tilde \phi_\theta (\eta)=\frac12 \eta^2 +O(\theta^2\eta^4).\] The function $\tilde \phi_\theta$  converges to $\tilde \phi_0$ in $C^N$ for any $N$ as $\theta\to 0$ on $\{ \eta:    |\eta| \lesssim 1    \}$.
Thus, if $\theta_0$ is small enough, the inequalities in Lemma \ref{wavepacket} hold uniformly in $\theta$ with $\mathcal T_\theta$ replacing $\mathcal R^*$  since the key geometric estimates are stable under smooth perturbation (in $C^N$ with a large enough $N$)  of the function $\tilde\phi_0$ (see \cite[Lemma 3.1]{L2}). 
Moreover, the $\epsilon$-removal argument also works in uniform manner since  
\[ |\mathcal T_\theta \psi(x,t)| \le C(1+|(x,t)|)^{-\frac{d-1}2}\] 
with a constant $C$, independent of $\theta$, for $\psi\in C_c^\infty((2^{-2}, 2^2)\times \mathbb R^{d-1})$.

\subsection{Deducing the linear estimate}  
To show  Proposition \ref{p-ls'}, we prove the following proposition,  making use of  Theorem \ref{bilinear-stronger}. Let us set 
\Be
\label{tBN}
\tilde B^\alpha_{p,q}(N,\theta)
= 
N^{2 (\frac{d-1}2+\frac{1}{p}- \frac{d}q)}
 \theta^{ 2((d-1)(1-\frac1p)- \frac{d-1+2\alpha}{q})}. 
 \Ee

\begin{prop}
\label{bi-1'}
Let $0< \alpha \le 1$ and  $ N^{-1/2} \le \theta \le \theta_0$ for a positive $\theta_0$ small enough.   Let $\mathcal E\subset [1,2]$ be $N^{-1}$-separated.  
 Suppose that  $ \Theta_1, \Theta_2$ are spherical caps such that    $\dist(\Theta_1, \Theta_2)\sim  \theta$ if 
 $\theta> 2 N^{-1/2}$ and $\dist(\Theta_1, \Theta_2)\le \theta$  if $\theta \le 2 N^{-1/2}$
  Suppose  that  $f_1$ and $f_2$  satisfy  \eqref{supp-con} and $\tilde{\mathcal A}_{\alpha}(\cE; 1/N)\le \tilde{\mathcal A}$ for a constant  $\tilde{\mathcal A}$.
  Then, there is a constant $C=C(\tilde{\mathcal A})$  such that 
\Be 
\label{bilineardiscrete'} 
\| e^{it\sqrt{-\Delta}} f_1 e^{it\sqrt{-\Delta}} f_2 \|_{L^{q /2}( \mathbb R^d\times \mathcal E)} \le C  \tilde B^\alpha_{p,q}(N,\theta) \|f_1\|_{p} \|f_2\|_{p}
\Ee
for $p,q$ satisfying 
\Be
\label{pq-con} 2/p=\bar q_\circ/q, 
\quad  \bar q_\circ> \tilde q_\circ.\Ee
\end{prop}

  Once we have the  estimate \eqref{bilineardiscrete'} for $p,q$ satisfying \eqref{pq-con}, 
 one can prove Proposition  \ref{p-ls'} by following the same argument in {\it Proof of  Proposition  \ref{p-ls}}.

\begin{proof}[Proof of Proposition  \ref{p-ls'}]
In order to show \eqref{localsmoothing-2d}, by finite decomposition of $f$ in the Fourier side, we may assume that  $\supp \widehat f$ is contained in  $\mathbb A_N$  and a sector with angular diameter smaller than 
$\theta_0$. So, we may apply Proposition \ref{bi-1'}. 
 
Indeed,  note that $q_\circ=2(d+3)/(d+1)> \tilde q_\circ$.  As before, let  $(p_\ast,q_\ast)$ be given by  $2/p_\ast=\bar q_\circ/q_\ast$ and 
\Be
\label{s-line*}
(d-1)\Big(1-\frac1{p_\ast} \Big)= \frac{d-1+2\alpha}{q_\ast}
\Ee
for some $\bar q_\circ \in ( \tilde q_\circ, q_\circ)$. 
Then, one can check that $1/p_\ast+1/q_\ast>1/2$ and $2< p_\ast< q_\ast$ hold.  
Therefore,  Proposition \ref{bi-1'} and Lemma  \ref{bisum} yield the estimate 
 \eqref{bnu} without  the factor $\mathcal A_{\alpha}^{2/q}$   for $(p,q)$ in a neighborhood of $(p_\ast, q_\ast)$ that satisfy  \eqref{pq-con}. 
Therefore, by  Lemma \ref{B-Sum-Tri} we obtain  
 \eqref{res-weak} without the factor  $\mathcal A_{\alpha}^{1/q_\ast}$ for $(p_\ast, q_\ast)$.

 Since we can take any $\bar q_\circ \in ( \tilde q_\circ, 2(d+3)/(d+1))$ for the estimate \eqref{res-weak} without the factor  $\mathcal A_{\alpha}^{1/q_\ast}$ obtained above, by real interpolation between this and the estimate \eqref{localsmoothing-2d} with $(p,q)=(1,\infty)$,  
we obtain the estimate  \eqref{localsmoothing-2d} for  $p,q$ satisfying  \eqref{s-line}  whenever $q>\tilde q_\ast$.
\end{proof}

Note that the estimate \eqref{bilineardiscrete'} is equivalent to \eqref{bilineardiscrete0} when $p=q=\infty$.  Thus, in order to prove Proposition \ref{bi-1'}, 
by interpolation    we only have to  show  
 \Be
 \label{bi-dis2} 
\| e^{it\sqrt{-\Delta}} f_1 e^{it\sqrt{-\Delta}} f_2 \|_{L^{q /2}( \mathbb R^d\times \mathcal E)} \lesssim  
\tilde B^\alpha_{2,q}(N,\theta)
 \|f_1\|_{2} \|f_2\|_{2}
 \Ee
 for  $q> \tilde q_\circ$  and $N^{-1/2}\le \theta \le \theta_0$ under the assumptions in Proposition \ref{bi-1'}.
 
\begin{proof}[Proof of  \eqref{bi-dis2}]    
 We only consider the case $\theta >  2 N^{-1/2}$. 
The case $\theta \le 2 N^{-1/2}$  can be shown by the same argument as in the proof of Proposition \ref{thm:bi}. Indeed, since  $\tilde{\mathcal A}_{\alpha}(\cE; 1/N)\le \tilde{\mathcal A}$,  we have $\#\mathcal{E}\le 
\tilde {\mathcal A} N^{\alpha}$. Since $\theta\sim N^{-1/2}$, it suffices to show 
\begin{equation}
\label{small-ang}
	\| e^{it\sqrt{-\Delta}} f_1 e^{it\sqrt{-\Delta}} f_2  \|_{L^{q/2}( \mathbb R^d)} \le C  N^{(d+1)(\frac1p-\frac{1}q)}   \|f_1\|_{L^p}\|f_2\|_{L^p}
\end{equation}
for all $t\in [1,2]$ 
whenever $ 2\leq p\leq q\leq \infty $. This estimate can be shown in the same manner as  \eqref{small-angle}, by interpolating the easy estimates for $(p,q)=(2,2), (2,\infty),$ and  $ (\infty,\infty)$.

 As before, by rotation we may assume that \eqref{sep} holds.    We now recall \eqref{ut} and \eqref{ut-lq}, and note that the estimate \eqref{bi-dis2} is equivalent to 
\Be
 \label{bilineardiscrete2'} 
\| U_t f_1 U_t f_2 \|_{L^{q /2}( \mathbb R^d\times \mathcal E)} \lesssim  
\tilde B^\alpha_{2,q}(N,\theta)
 \|f_1\|_{2} \|f_2\|_{2}.
 \Ee
 To show this, we use Corollary  \ref{bilinear-theta} instead of Theorem \ref{thm:cone}.  

Thanks to the support condition  \eqref{sep}, we may  use the  locally constant property at scale $N^{-1}\theta^{-2 }$.  By \eqref{lc-theta}, we have 
\Be 
\label{ineq-}   \Big( \sum_{ t\in \cE}  \| U_t f_1 U_t f_2 \|_{q/2}^{q/2} \Big)^{2/q} \lesssim   \sigma^\frac2q  \Big( \sum_{ t\in 
\cE( N^{-1}\theta^{-2 })}  \| \mathfrak S(U_t f_1) \mathfrak S(U_t f_2) \|_{q/2}^{q/2} \Big)^{2/q},    \Ee
where $\sigma$ is given by  \eqref{equ-sig}. Here, as before, $\cE(\delta)$ denotes a maximally $\delta$-separated subset of $\cE$.  As it becomes clear later, the estimate 
holds uniformly  regardless of the choice of the maximally separated set. 
 
We set 
\[ 
U_t^\theta h= \int e^{i(x\cdot\xi+ t\xi_1\tilde \phi_\theta(\xi'/\xi_1) )} \widehat h(\xi) \,d\xi, 
\]
where $\tilde \phi_\theta$ is given by \eqref{phi-theta}. By the standard parabolic rescaling, we have 
\[ (U_t f_1 U_t f_2)(x)=(U_{\theta^2 t}^\theta f_{1,\theta}  U_{\theta^2t}^\theta f_{2,\theta})(x_1,\theta x'), \] 
where $\widehat {f_{j,\theta}}(\xi)=\theta^{d-1} \widehat f_j(\xi_1,\theta \xi')$, $j=1,2$.  This observation and changing variables $(x',t)\to (\theta^{-1} x',\theta^{-2} t)$ for the right hand side 
of \eqref{ineq-}  give
\[    \Big( \sum_{ t\in \cE}  \| U_t f U_t g \|_{q/2}^{q/2} \Big)^{2/q} \lesssim   \theta^{-\frac{2(d-1)}q} \sigma^\frac2q  \Big( \sum_{ t\in \theta^2 \cE( N^{-1}\theta^{-2 })}  \| \mathfrak S(U_t^\theta f_{1,\theta}) \mathfrak S(U_t^\theta f_{2,\theta})  \|_{q/2}^{q/2} \Big)^{2/q}.  \] 

Let us set 
\[  \mathfrak T= N\theta^2 \cE( N^{-1}\theta^{-2 }).\]
After rescaling  $(x,t)\to  N^{-1}(x,t)$ and  $\xi\to N\xi$ in the frequency side  for both $U_t^\theta f_{1,\theta}$ and $U_t^\theta f_{2,\theta}$, we bound the right hand side of the above inequality  by 
\[    N^{2(\frac{d}2 -\frac dq)}  \theta^{d-1-\frac{2(d-1)}q} \sigma^\frac2q  \Big( \sum_{ t\in \mathfrak T}\|\mathfrak S(U_t^\theta  f_{1,\theta}') \mathfrak S(U_t^\theta  f_{2,\theta}')  \|_{q/2}^{q/2} \Big)^{2/q} \] 
while $\|f_{j,\theta}'\|_2=\|f_j\|_2$ for $j=1,2$. In fact, $\widehat {f_{j,\theta}'}=  \theta^{-\frac{d-1}2} N^{\frac d2}\widehat{f_{j,\theta}}(N\cdot)$. 
Let us set 
\[ \tilde f_{j,\theta}= \widehat {f_{j,\theta}'}, \quad j=1,2.\] 
Consequently,  we have 
 $\| \tilde f_{1,\theta}\|_2=\|f_1\|_2$,  $\| \tilde f_{2,\theta}\|_2=\|f_2\|_2$ by Plancherel's theorem. Also, it follows that
 \Be
 \label{sep-1}
 \supp \tilde f_{1,\theta}  \subset \{ \xi: \xi_1\sim 1,  |\xi'| \sim 1      \},      		\quad  		\supp \tilde f_{2,\theta}  \subset \{ \xi_1\sim 1,  |\xi'| \ll 1\}.
 \Ee

From \eqref{R*-theta}, 
note that $U_t^\theta  f_{j,\theta}'= \mathcal T_\theta\tilde f_{j,\theta}$, $j=1,2.$ 
Recalling \eqref{sig} and \eqref{tBN}, we  obtain 
\begin{align*}
   \Big( \sum_{ t\in \cE}  \| U_t f_1U_t f_2 \|_{q/2}^{q/2} & \Big)^{2/q} \lesssim   \tilde B_{2,q}^\alpha(N,\theta ) 
  \Big( \sum_{ t\in \mathfrak T}  \| \mathfrak G(\mathcal T_\theta \tilde f_{1,\theta}(\cdot ,t))  \mathfrak G(\mathcal T_\theta\tilde f_{2,\theta}(\cdot ,t))  \|_{q/2}^{q/2} \Big)^{2/q}. 
   \end{align*}

Note that the set $\mathfrak T=N\theta^2 \cE( N^{-1}\theta^{-2 })$ is 1-separated. 
Let  $\mathbf I_t=[t, t+1]$ for $t\in \mathfrak T$ and consider the union 
\[ \Gamma=\bigcup_{t\in \mathfrak T} \, \mathbf I_t .\]  \
Since $\tilde f_{1,\theta}$ and $\tilde f_{2,\theta}$ are compactly supported,  both $ \mathcal T_\theta\tilde f_{1,\theta}(\cdot ,t)$ and $\mathcal T_\theta\tilde f_{2,\theta}(\cdot ,t)$ enjoy the  locally constant property in $t$ at scale $1$. 
Indeed,  for $|t-t'|\le 1$, $\xi_1\sim 1$, and $|\xi'|\lesssim 1$,  we note $|\partial_\xi^\alpha((t-t')  \xi_1\tilde \phi_\theta(\xi'/\xi_1))|\lesssim_\alpha 1 $
(see the proof of Lemma~\ref{lcp}). Hence, for any $t\in \mathfrak T$, we have 
\[   
\| \mathfrak G(\mathcal T_\theta\tilde f_{1,\theta}(\cdot,t))  \mathfrak G(\mathcal T_\theta\tilde f_{2,\theta}(\cdot,t))  \|_{L^{q/2}_x}\lesssim \| \mathfrak G^2(\mathcal T_\theta\tilde f_{1,\theta}(\cdot ,s))  \mathfrak G^2(\mathcal T_\theta\tilde f_{2,\theta}(\cdot,s))  \|_{L^{q /2}_{x,s}( \mathbb R^d\times \mathbf I_t)}.  \]
Here, as mentioned before,  we are actually abusing the notation $\mathfrak G$ since $\mathfrak G^2$ is not  a composition of the same operator $\mathfrak G$. Consequently, combining this and the preceding  inequality, we have
\begin{align*}
   \Big( \sum_{ t\in \cE}  \| U_t f_1U_t f_2 \|_{q/2}^{q/2}  \Big)^{2/q} \lesssim   \tilde B_{2,q}^\alpha(N,\theta )  
   \Big\| \mathfrak G^2(\mathcal T_\theta\tilde f_{1,\theta}(\cdot ,s))  \mathfrak G^2(\mathcal T_\theta\tilde f_{2,\theta}(\cdot ,s)) \Big \|_{L^{q /2}_{x,s}( \mathbb R^d\times \Gamma)}.
   \end{align*}
   
Therefore, to complete the proof of  \eqref{bilineardiscrete2'}, we only have to show 
\Be \label{fin} \Big\|  \mathcal T_\theta \tilde f_{1,\theta}  \mathcal T_\theta \tilde f_{2,\theta} \Big \|_{L^{q /2}( \mathbb R^d\times \Gamma)}
\le C(\tilde{\mathcal A}) \| \tilde f_{1,\theta}\|_2\|\tilde f_{2,\theta}\|_2
\Ee
for $q> \tilde q_\circ$  and $N^{-1/2}\le \theta \le \theta_0$ when \eqref{sep-1} holds.  To show this, we use the uniform estimate in Corollary  \ref{bilinear-theta}. For this purpose, it is sufficient to show that the collection $\mathcal I$ satisfies \eqref{intervals}. Furthermore, since the set $\mathfrak T$ is separated by $1$,  it is enough to show 
\Be \label{aa} \#(\mathfrak T \cap I_\rho) \le  \tilde {\mathcal A} \rho^\alpha\Ee
for  $N^{-1/2} \le \theta \le \theta_0$  whenever  $I_\rho$ is an interval of length $\rho \geq 1$.  Indeed, once we have this, the estimate \eqref{fin} follows from Corollary  \ref{bilinear-theta}. 

Recalling $\mathfrak T=N\theta^2 \cE( N^{-1}\theta^{-2 })$ and $\cE( N^{-1}\theta^{-2 })\subset I_\circ=: [1,2]$, we note that 
\begin{align*}
\#(\mathfrak T \cap I_\rho)
&=\#(\cE( N^{-1}\theta^{-2 })\cap (N^{-1}\theta^{-2 } I_\rho\cap I_\circ)) 
\\
&=  \Big( \frac{N^{-1}\theta^{-2 }}{|I'|} \Big)^\alpha \mathcal N (\cE\cap I',  N^{-1}\theta^{-2 })  \Big( \frac{|I'|}{N^{-1}\theta^{-2 }}\Big)^\alpha 
\end{align*}
where $I'=N^{-1}\theta^{-2 }I_\rho \cap I_\circ $. Consequently,  from  the definition  \eqref{agamma'} we see that
\[   \#(\mathfrak T \cap I_\rho)\le   {\mathcal A}_\alpha(\cE; N^{-1}\theta^{-2 })   \rho^\alpha \le   \tilde{\mathcal A}_\alpha(\cE; 1/N)   \rho^\alpha.   \]
Therefore, \eqref{aa} follows. 
\end{proof}

\begin{rem}  If $\theta\sim 1$ in Proposition \ref{bi-1'},  the estimate \eqref{bilineardiscrete'} holds   with $C=C(\mathcal A)$ 
under a weaker assumption that 
${\mathcal A}_{\alpha}(\cE; 1/N)\le  {\mathcal A}$ for a constant  ${\mathcal A}$. It is clear from the last part of the proof. 
\end{rem} 

\subsection{Sharp estimates with $\epsilon$-loss} 
\label{sec:sss}
Until now, we have been concerned with the estimate \eqref{smooth-pq} with the optimal exponent $ s = s_c(p,q) $.  
However, if an $ \epsilon $-loss is allowed, it is possible to obtain estimates over an extended range of exponents.  
To keep the discussion focused, we restrict our attention to the estimate
\Be \label{smooth-pq0}
\left( \sum_{t \in E(2^{-j})} \left\| e^{it\sqrt{-\Delta}} P_j f \right\|_{L_x^q(\mathbb{R}^d)}^q \right)^{\frac{1}{q}} 
\le C\, 2^{(s_c(p,q) + \epsilon) j} \|f\|_p,
\Ee
for $ p $ and $ q $ satisfying \eqref{s-line}, which, as mentioned earlier, is the most interesting case.

Let us consider the estimate 
\Be 
\label{localsmoothing0} 
\| e^{it\sqrt{-\Delta}} f\|_{L^r( \mathbb R^d\times I_\circ)} \le C N^{\frac{d-1}2-\frac{d}{r}+\epsilon} \|f\|_r 
\Ee
for any $\epsilon>0$ while $\supp \widehat f\subset \mathbb A_N$. It has been conjectured that the estimate holds for $r\ge 2d/(d-1)$. 
When $d=2$, this is verified by Guth--Wang--Zhang \cite{GWZ}. In higher dimensions, the estimate \eqref{localsmoothing0} was shown by Bourgain--Demeter for $r\ge 2(d+1)/(d-1)$. 
Recently, the range  is further extended  by Gan--Wu \cite{GW}.

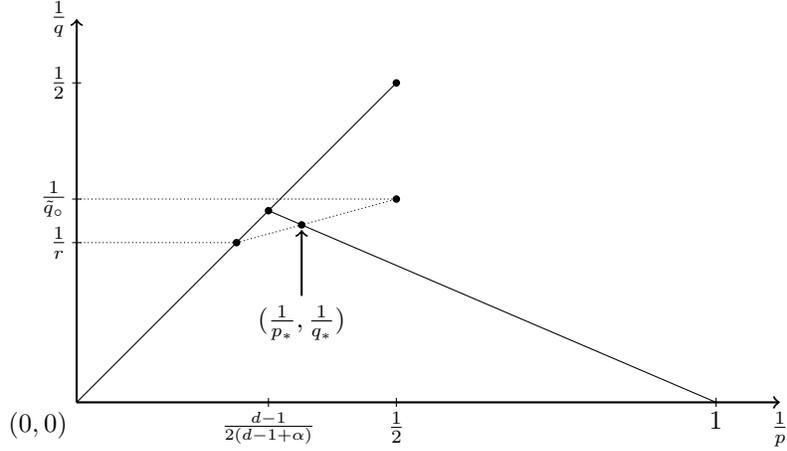
\begin{figure}[t] 
	\centering 
	\begin{tikzpicture}[scale=0.85]
		
		\draw[thick,->]  (0,0) -- (11,0);     
		\draw[thick,->] (0,0) -- (0,6);
		\node[below] at (11,0) {$\frac{1}{p}$};
		\node[left] at (0,6) {$\frac{1}{q}$};
		
		\node[below left] at (0,0) {$(0,0)$};  
		
		\draw (10,-2pt) -- (10,2pt);           
		\node[below] at (10,0) {$1$};          
		
		\draw (-2pt,5) -- (2pt,5);  
		\node[below] at (5,0) {$\frac{1}{2}$};
		
		\draw (5,-2pt) -- (5,2pt);             
		\node[left] at (0,5) {$\frac{1}{2}$}; 
		
		\draw (-2pt,2.5) -- (2pt,2.5);             
		\node[left] at (0,2.5) {$\frac{1}{r}$}; 
		
		\draw (-2pt,35/11) -- (2pt,35/11);        
		\node[left] at (0,35/11) {$\frac{1}{\tilde q_\circ}$}; 
		
		\draw (3,-2pt) -- (3,2pt); 
		\node[below,font=\footnotesize] at (3,0) { $\frac{d-1}{2(d-1+\alpha)}$};

		\draw[->, thick] (95/27,15/9) -- (95/27,24.3/9);   
		\node[below] at (95/27,15/9) {$(\frac{1}{p_\ast},\frac{1}{q_\ast})$}; 
		
		
		

		\node[circle, fill=black, inner sep=1.0pt] at (5,35/11) {};
		\node[circle, fill=black, inner sep=1.0pt] at (5,5) {};
		\node[circle, fill=black, inner sep=1.0pt] at (95/27,25/9) {};
		\node[circle, fill=black, inner sep=1.0pt] at (3,3) {};
		\node[circle, fill=black, inner sep=1.0pt] at (2.5,2.5) {};

		\draw[densely dotted,black] (0,35/11) -- (5,35/11); 
		\draw[densely dotted,black] (0,2.5) -- (2.5,2.5); 
		\draw (0,0) -- (5,5);     
		\draw[densely dotted,black] (5/2,5/2) -- (5,35/11);      
		\draw (3,3) -- (10,0);
	\end{tikzpicture}
	\caption{ Using $L^r$  estimate \eqref{locals-} in place of the $L^\infty$ estimate  yields the sharp (with $\epsilon$-loss) estimate  on an extended range. } 
	\label{Fig3} 
\end{figure}

The estimate  can be used to extend the  $p,q$ range for \eqref{smooth-pq0} (see Figure \ref{Fig3}). Let $ \cE $ be a $1/N$-separated subset of $[1,2]$.  Then, it is easy to see that  
\Be 
\label{locals-} 
\| e^{  it\sqrt{-\Delta}} f\|_{L^r( \mathbb R^d\times  \cE)} \le C N^{\frac{d-1}2-\frac{d-1}{r}+\epsilon} \|f\|_r
\Ee  
for $\supp \widehat f\subset \mathbb A_N$.  This follows from the locally constant property at  scale $N^{-1}$. Note that 
$|e^{  is\sqrt{-\Delta}} f(x)|\le \mathfrak S( e^{  it\sqrt{-\Delta}} f)(x)$ if $|t-s|\le 1/N$. 
 Thus,  for any interval $I$ of length $1/N$ contain $s$, we have $\| e^{is\sqrt{-\Delta}} f\|_{L^r( \mathbb R^d)}\le CN^{1/r}\| e^{  it\sqrt{-\Delta}} f\|_{L^r( \mathbb R^d\times I)}$.  Since 
$ \cE $ is a $1/N$-separated subset of $[1,2]$,  it follows that 
\[    \| e^{  it\sqrt{-\Delta}} f\|_{L^r( \mathbb R^d\times  \cE)} \le C N^{1/r} \| e^{  it\sqrt{-\Delta}} f\|_{L^r( \mathbb R^d\times I_\circ)}.    \]
Therefore, \eqref{locals-} follows from \eqref{localsmoothing0}.

Using the estimate \eqref{locals-} together with the bilinear estimate \eqref{bilineardiscrete'} from Proposition~\ref{bi-1'} (with $ p = 2 $ and $ \theta \sim 1 $), one can, by following the argument in the proof of Proposition~\ref{bi-1'}, obtain the estimate
\Be
\label{bi-e}
\left\| \prod_{j=1}^2 e^{it\sqrt{-\Delta}} f_j \right\|_{L^{q/2}(\mathbb{R}^d \times \mathcal{E})}
\le C N^{\epsilon} \tilde{B}^\alpha_{p,q}(N, \theta) \prod_{j=1}^2 \|f_j\|_p,
\Ee
for all $ p, q $ satisfying
\Be
\label{pq-con0}
2 \le p \le r, \quad
\frac{2(r - \tilde{q}_\circ)}{ \tilde{q}_\circ(r - 2)} \left( \frac{1}{2} - \frac{1}{p} \right) < \frac{1}{\tilde{q}_\circ} - \frac{1}{q},
\Ee
whenever $ f_1 $, $ f_2 $, and $\cE$ satisfy  the assumptions  in Proposition~\ref{bi-1'}. 
In fact, due to the presence of $\epsilon$-loss in the exponent of $N$,  there is actually a $c\epsilon$-gain in the exponent of  $\theta$.

Indeed, by H\"older’s inequality, the estimate \eqref{locals-} may be viewed as a bilinear estimate \eqref{bi-e} with $ p = q = r $ and $ \theta \sim 1 $. Similarly, the estimate \eqref{bi-e} with $ q = \infty $ and $ \theta \sim 1 $ holds for all $ 1 \le p \le \infty $, as already  established in Propositions \ref{p-ls} and \ref{bi-1}.    Interpolation between those estimates  and \eqref{bilineardiscrete'} with $p=2$ and $\theta\sim 1$ yields  \eqref{bi-e} with $\theta\sim 1$
for $p,q$ satisfying \eqref{pq-con0} whenever 
$f_1$ and $f_2$  meet  the same condition as in Proposition \ref{bi-1'} with $\theta\sim 1$ and ${\mathcal A}_{\alpha}(\cE; 1/N)\le {\mathcal A}$ for a constant  $\mathcal A$. 
Then, the argument used in the proof of Proposition \ref{bi-1'} can be repeated to obtain \eqref{bi-e}.

Now, using the estimate \eqref{bi-e}\footnote{To be strict, we need bilinear estimates for $U_t^\theta$ in place of $e^{it\sqrt{-\Delta}}$
with bound uniform in $\theta$. However, since $(\xi',  \xi_1\tilde \phi_\theta(\xi'/\xi_1) )$ is a small smooth perturbation  of the standard cone, all the known local smoothing results for $e^{it\sqrt{-\Delta}}$ continue to be valid for  $U_t^\theta$ with bounds uniform in $\theta$.} and Lemma \ref{bisum} \footnote{Thanks to the currently known range of the sharp local smoothing estimate \eqref{locals-} (see, e.g., \cite{BD, GWZ, GW}), we may assume $r\leq 4$. As a result, the pair $(p_\ast, q_\ast)$ determined below satisfies $1/2 < 1/p_\ast + 1/q_\ast$ and $2 < p_\ast < q_\ast$ for all $0 < \alpha \leq 1$ and $d \geq 2$, so we may apply Lemma~\ref{bisum}.} along with the decomposition \eqref{bi-decom}, one can obtain \eqref{smooth-pq0} for $p, q$ satisfying \eqref{s-line} and  $q> q_\ast(\alpha, r)$, where $q_\ast:=q_\ast(\alpha, r)$ is determined by solving   
\eqref{s-line*} and 
\[ \frac{ 2( r- \tilde q_\circ) }{ \tilde q_\circ(r-2)} \Big(\frac12-\frac1p_\ast\Big)=  \frac1{\tilde q_\circ}  -\frac1q_\ast.\] 
(See Figure \ref{Fig3}.)  This can be done by reiterating the proof of Proposition \ref{p-ls} (also see  that of Proposition  \ref{p-ls'}). Since the desired estimates  for $(f,g)\mapsto B^{\nu_\circ} (f,g)$  follow  from \eqref{small-ang},  one  needs only to consider  the operators 
$\{B^{\nu_\circ}\}_{0\le \nu<\nu_\circ}$.

Indeed, a computation shows 
\[ 
q_\ast(\alpha, r)=\frac{r\big(2(d-1+2\alpha)^2-4\alpha^2\big)-4\big(4\alpha^2+(d-1)^2+5\alpha(d-1)\big)}{r(d-1)(2\alpha+d-1)-2(d-1)(3\alpha+d-1)}.
\]
 In particular, when $d=2$, by the sharp $L^4$ local smoothing estimate \cite{GWZ},   we may take $r=4$. 
 Consequently, we  have 
\eqref{smooth-pq0}  for $p, q$ satisfying \eqref{s-line} and 
\[ q> \frac{2+6\alpha}{1+\alpha} .\]
Moreover, taking  $r\to \infty$,  we note that the range $q>q_\ast(\alpha, \infty)$ matches the one  in Theorem \ref{improved}.

\begin{rem}\label{hns}
Finally, for $d\ge 4$,  we remark that  the $\epsilon$-loss in \eqref{smooth-pq0} can be removed in a certain range of $p,q$, which is even broader  than the range in Theorem \ref{improved}.  The optimal local smoothing estimate \eqref{localsmoothing0} without the loss $N^{\epsilon}$  was  obtained by 
Heo--Nazarov--Seeger \cite{HNS} for $r>  r_d:=2+\frac4{d-3}$ when $d\ge 4$. Applying this result, one can eliminate  the $2^{\epsilon j}$-loss in \eqref{smooth-pq0}. Indeed, we obtain the estimate \eqref{smooth-pq} with $s=s_c(p,q)$  for $p,q$ satisfying 
\eqref{s-line} and $q>q_\ast(\alpha, r_d)$. 
\end{rem}

\section{Sharpness of the regularity exponent} \label{sec:lower}	
In this section,  we discuss the 
sharpness of the regularity exponent $s$ in the estimate \eqref{smooth-pq}. 

	We begin  by recalling the definitions \eqref{agamma}  and \eqref{agamma'} for a discrete set $\mathcal E$. The following lemma shows the existence of a certain $2^{-j}$-separated set,  which will  be used in the proof of the lower bounds in Proposition \ref{prop:neces}.

\begin{lem}
\label{lem:existA}
Let $ 0\leq \alpha\leq 1 $. For each $j\ge 1$, there exists a  $2^{-j}$-separated subset $\mathcal{E}\subset [1,2] $  such that $\#\mathcal{E}\sim 2^{j\alpha}$, $ \mathcal A_{\alpha}(\mathcal E;  2^{-j})\lesssim 1$,   and $ \tilde{\mathcal A}_{\alpha}(\mathcal E;  2^{-j})\lesssim 1$.
\end{lem}

\begin{proof}
If $ \alpha=0 $, we take $ \mathcal{E}$ to be  a set of one element. Thus it remains to consider the case $ 0<\alpha\leq 1 $. 

Let $ \mu=2^{-1/\alpha}\in (0,1/2]$. We use the standard construction for a Cantor set, see \cite[p. 60]{1995mattila}. At step $ 0 $, set $ I_{0,1}^\mu=[1,2]$. At step $ 1 $, we take a closed interval $ I_{1,1}^\mu$ of length $ \mu $ such that it shares the same left endpoint with the interval $ I_{0,1}^\mu$. Similarly, we let $ I_{1,2}^\mu$ be the interval of length $ \mu $ such that it shares the same right endpoint with the interval $ I_{0,1}^\mu$.  We continue this process for the  closed intervals appearing in the previous step. At stage $ k $, we have $ 2^k $ intervals $ I_{k,1}^\mu,\cdots, I_{k,2^k}^\mu$ of length $ \mu^k $. 

We now choose $k$ such that  $ \mu^k\sim L2^{-j} $ with $ L $ being a fixed large positive constant.  We let $ \mathcal{E} $ be the set of all the right endpoints of these $ 2^k $ intervals $ I_{k,1}^\mu,\cdots, I_{k,2^k}^\mu$. Clearly, $ \mathcal{E} $ is $2^{-j}$-separated with cardinal $\#\mathcal{E}\sim 2^{j\alpha} $. 

To see $ \mathcal A_{\alpha}(\mathcal E;  2^{-j})\lesssim 1$,  let $I\subset [1,2]$ be an interval of length $\sim \mu^l\ge 2^{-j}.$ From the construction, it is clear that  $\mathcal N(\mathcal E\cap I, 2^{-j})$ maximizes when $I$ shares its left or right endpoint with $[1,2]$. In such cases, we have $\mathcal N(\mathcal E\cap I, 2^{-j})\sim 2^{k-l}$. Therefore, 
\[  \mathcal A_{\alpha}(\mathcal E;  2^{-j})\lesssim \Big( \frac{\mu^k}{\mu^l}\Big)^\alpha  2^{k-l}\lesssim 1. \]
Similarly, one can easily see from the self-similar nature of the construction  that ${\mathcal A}_{\alpha}(\mathcal E;  2^{-m})\lesssim 1$ for all $1\le m\le j$.  Hence, 
it follows that $\tilde {\mathcal A}_{\alpha}(\mathcal E;  2^{-j})\lesssim 1$.
\end{proof}

\begin{rem}[Explicit   description of the construction in the proof of Lemma \ref{lem:existA}]
\label{rem:description}  It is not hard to see that, for each $k\in \mathbb{N}$,
\[\mathcal{E}=\Big\{1+\mu^k+\sum_{m=0}^{k-1}(1-\mu)\mu^m P_m\, :\, P_m=0,1\Big\}.\]
This gives a decomposition  $ \mathcal{E} =\bigcup_{l=0}^k\mathcal{E}_l$, where  $\mathcal{E}_k=\{1+\mu^k\}$  and 
\[\mathcal{E}_l=\Big\{1+\mu^k+(1-\mu)\mu^{l}+\sum_{m=l+1}^{k-1}(1-\mu)\mu^m P_m\, :\, P_m=0,1\Big\}, \quad 0\leq l\leq k-1.\]
Then, it follows that $\#\mathcal{E}_l=2^{k-l-1}$ for  $ 0\leq l \leq k-1 $. Observe also that $ \mu^k\sim L2^{-j} $ implies $ k=j\alpha+O(1) $. 
\end{rem}

Before we show  the 
sharpness of the regularity exponent $s$ in the estimate \eqref{smooth-pq}, we first consider the case of discrete dilation sets. 
The desired conclusion will easily  follow from this result for the discrete setting.

	\begin{prop}
	\label{prop:neces}
		Let $ 0\leq \alpha\leq 1$, $ 1\leq p\leq q\leq \infty$, and  $\mathcal E\subset [1,2]$ be a $2^{-j}$-separated set. Suppose that  the estimate
	\Be 	
	\label{hohoho}
			\|e^{it\sqrt{-\Delta}}P_jf\|_{L^q(\mathbb{R}^d\times \mathcal{E})}\le C2^{js}\|f\|_{L^p(\mathbb{R}^d)}
	\Ee		
holds for a constant $C$, independent of $j\ge 1$, whenever $\mathcal A_{\alpha}(\mathcal E;  2^{-j})\le  C_0 $ for a constant $C_0>0$. Then,  we must have
		\[
			s\geq  s_c(p,q):=\max  \big\{s_1(p,q),\,s_2(p,q),\, s_3(p,q)\big\}. 
	\]
		Moreover, if $ p=q=\frac{2(\alpha+d-1)}{d-1}$ and $\alpha>0$, then it is necessary that $s>{\alpha}/{q}$.
	\end{prop}
	
	The examples used in the proof of this proposition are adaptations of  the standard examples (see, e.g., \cite[Section 4]{TV2}).
	
	\begin{proof}
	We prove each of the lower bounds  $s\ge s_1(p,q)$,  $s\ge s_2(p,q)$, and $s\ge s_3(p,q)$,   separately.
	Without loss of generality,  we can replace $ P_jf $ with $ f $ whose Fourier support is contained in $\mathbb{A}_{2^j}=\{\xi: 2^{j-1}\leq|\xi|\leq 2^{j+1}\}$.

\noindent \textit{Lower bound $s\geq  s_1(p,q):= \frac{d-1}{2}+\frac{1}{p}-\frac{d}{q}$.}
	As  mentioned in the introduction,  this lower bound applies to  a fixed time estimate  \eqref{single} for  the wave operator. 
	 
	 Let $\mathcal{E}=\{1+L2^{-j}\}$ with $j\gg 1$, which trivially satisfies $\ \mathcal A_{\alpha}(\mathcal E;  2^{-j})\le 1$.  Here $ L $ is a fixed large positive constant as in the proof of Lemma \ref{lem:existA}. We take $f$  given by 
	 \Be\label{hhh} \widehat{f}(\xi)=e^{-i|\xi|}\beta_1(|\xi|/2^j) ,\Ee where $ \beta_1 \in C_c^\infty\big((1/4,4)\big) $ is a non-negative  function such that $ \beta_1=1 $ on $ [2^{-1},2]$. Since $\widehat{f}$ is radial, the Fourier inversion formula for radial functions \cite[p.~155]{SW} gives 
		\[f(x)=c2^{j\frac{d+2}{2}}\int_{0}^{\infty}e^{-i2^jr}\beta_1(r)J_{\frac{d-2}{2}}(2^jr|x|)(r|x|)^{-\frac{d-2}{2}}r^{d-1}\,dr\]
		for a constant $c$,  where $ J_{\frac{d-2}{2}}$ is the Bessel function. It is well known that 
		\[J_{\frac{d-2}{2}}(r)=e^{\pm ir}a_{\pm}(r),\]
		where $ a_\pm\in S^{{-1/2}} $ in the sense that $|(\tfrac{d}{dr})^ma_\pm(r)|\lesssim (1+r)^{-m-{1/2}}$ for any $ m\in\mathbb{N}$ (see \cite[Chapter VIII, 1.4.1]{stein1993harmonic}). This gives 
		\begin{align*}
			f(x)&=c2^{j\frac{d+2}{2}}\sum_\pm \int_{0}^\infty e^{i2^jr(\pm |x|-1)}\beta_1(r)a_\pm(2^jr|x|)(r|x|)^{-\frac{d-2}{2}}r^{d-1}\,dr.
		\end{align*}
		By routine integration by parts, it follows that  $|f(x)|\lesssim 2^{j\frac{d+1}{2}}(1+2^j||x|-1|)^{-N}$
		for any $ N\in\mathbb{N}$.  Consequently, it follows that   
		\[  \|f\|_{L^p(\mathbb{R}^d)} \lesssim 2^{j(\frac{d+1}{2}-\frac{1}{p})}.\]
		
		To obtain a lower bound for $e^{it\sqrt{-\Delta}}f(x) $, we use  the asymptotic  expansion 
		\[J_\frac{d-2}{2}(r)=c_\pm e^{\pm ir}r^{-1/2}+R(r),  \quad r\geq 1,\]
		where $ c_{\pm}\neq 0 $ and $ |R(r)|\lesssim r^{-3/2}$ (see \cite[Chapter VIII, 5.2]{stein1993harmonic}). 
		
		Let us set 
		\begin{align*} 
		I_\pm(x,t) &=  c_\pm 2^{j\frac{d+1}{2}}|x|^{-\frac{d-1}{2}}\int e^{i2^jr(t-1\pm |x|)}\beta_1(r)r^{\frac{d-1}{2}}\,dr, 
		\\ 
		I_0(x,t)&=2^{j\frac{d+2}{2}}|x|^{-\frac{d-2}{2}}\int e^{i2^jr(t-1)}\beta_1(r)R(2^jr|x|)r^{\frac{d}{2}}\,dr.
		\end{align*}
		Thus, we have
		\begin{align*}
			e^{it\sqrt{-\Delta}}f(x)&=  I_+(x,t) +I_-(x,t)+  I_0(x,t) 	               
		\end{align*}
		when $ |x|\gtrsim 2^{-j} $.  It is easy to see that  $I_0(x,t) = O(2^{j\frac{d+1}{2}}|x|^{-\frac{d-1}{2}} (2^j|x|)^{-1})$ for $ |x|\gtrsim 2^{-j} $.  
		Also note that  
		\[ I_\pm(x,t) =c_\pm 2^{j\frac{d+1}{2}}|x|^{-\frac{d-1}{2}}\phi\big(2^j(t-1\pm |x|)\big),\]
		where   $ \widehat{\phi}(s)=\beta_1(s) s^{(d-1)/2} $.  
		 We take $t=1+L2^{-j} \in\mathcal{E}$ and let 
		 \[ J_t= [t-1-c_02^{-j},t-1].\] 
		Noting that $ \phi\in \mathcal{S}(\mathbb{R}) $ and $ \phi(0)>0 $,  we  choose $ c_0 $ small enough and $ L $ large enough so that $ |I_{+}(x,t)|\ll |I_{-}(x,t)|\sim 2^{j\frac{d+1}{2}}|x|^{-\frac{d-1}{2}}$ if $|x|\in J_t$. Thus, we obtain 
		 \Be
		 \label{hh} |e^{it\sqrt{-\Delta}}f(x)|\gtrsim2^{j\frac{d+1}{2}}|x|^{-\frac{d-1}{2}} \Ee
		 provided that $|x| \in J_t$. As a result, we have
		\begin{align*}
		\|e^{it\sqrt{-\Delta}}f\|_{L^q(\mathbb{R}^d)}\gtrsim 2^{j\frac{d+1}{2}}\Big(\int_{J_t}r^{-\frac{d-1}{2}q}r^{d-1}\,dr\Big)^{1/q}\gtrsim 2^{jd}2^{-j\frac{d}{q}}.
		\end{align*}
		Since $  \|f\|_{L^p(\mathbb{R}^d)} \lesssim 2^{j(\frac{d+1}{2}-\frac{1}{p})}$, the estimate  \eqref{hohoho} with $\mathcal{E}=\{1+L2^{-j}\}$ implies 
		\[2^{jd}2^{-j\frac{d}{q}}\le C 2^{j(s+\frac{d+1}{2}-\frac{1}{p})},\]
		which yields the desired bound $s\geq  s_1(p,q)$  by  letting $ j\to\infty $.

\noindent \textit{Lower bound $s\geq s_2(p,q):= \frac{d+1}{2}\big(\frac{1}{p}-\frac{1}{q}\big)+\frac{\alpha}{q}$.}
Using Lemma \ref{lem:existA}, we can take a $2^{-j}$-separated set $ \mathcal{E}\subset [1,2] $  with cardinality  $\#\mathcal{E}\sim 2^{j\alpha} $ and $ \mathcal A_{\alpha}(\mathcal E;2^{-j})\lesssim 1$.  For each $t\in \mathcal E$, we set 
\[ S_t=\{x:  |x'|\le c 2^{-j/2}, \quad  |x_d+t|\le c 2^{-j}\} \]
with  a small positive constant $c$, which we choose later. 

Recall that $\beta_0\in C_c^\infty((-2^2, 2^2))$ with  $\beta_0=1$ on $[-2, 2]$,  and  $\beta_1 \in C_c^{\infty}((2^{-2}, 2^2))$  with  $\beta_1=1$ on $[2^{-1}, 2]$. 
 Let us consider $f$ given by 
\[\widehat{f}(\xi)=\beta_0\Big(\frac{|\xi'|}{c_12^{j/2}}\Big) \beta_1\Big(\frac{\xi_d}{2^j}\Big), \quad \xi=(\xi',\xi_d), \] which is essentially the Knapp example. Here $ c_1 $ is a small positive constant to be chosen later. It is easy to see that $ \|f\|_{L^p(\mathbb{R}^d)} \sim 2^{j\frac{d+1}{2}(1-\frac1p)}$.  Scaling  $(\xi^{\prime}, \xi_d)\to (2^{j/2}\xi', 2^j  \xi_d)$ gives 
\[e^{it\sqrt{-\Delta}}f(x)=\frac{2^{j(d+1)/2}}{ (2\pi)^{d} }\!\int e^{i2^{j/2}x'\cdot\xi'}e^{i2^j(x_d+t)\xi_d}e^{it\xi_d\tilde \phi_{2^{-j/2}}(\xi'/\xi_d)}\beta_0\Big(\frac{|\xi'|}{c_1}\Big) \beta_1 ({\xi_d})\,d\xi' d\xi_d,\]
where $\tilde \phi_{2^{-j/2}}$ is given by  \eqref{phi-theta}. Note  that $|e^{it\xi_d\tilde \phi_{2^{-j/2}}(\xi'/\xi_d)}-1|\lesssim  2|\xi'|^2/{\xi_d}\lesssim c_1^2$ 
if $ t\in[1,2] $ and   $\xi\in \supp \widehat f(2^{j/2}\cdot,2^{j}\cdot)$. Thus, if $c$  and $ c_1 $ are small enough, we have 
\[ |e^{it\sqrt{-\Delta}}f(x)|\gtrsim 2^{j\frac{d+1}2}\]
 for any $ t\in [1,2]$  if $x\in  S_t$. 
 Hence, $\|e^{it\sqrt{-\Delta}}f\|_{L^q(\mathbb{R}^d)}\gtrsim 2^{j(\frac{d+1}2 -\frac{d+1}{2q})}$ for $ t\in [1,2]$. 
 
 Since  $\# \mathcal E\sim 2^{\alpha j}$,
 the left hand side of the inequality \eqref{hohoho} is bounded below by  $2^{j(\frac{d+1}2 -\frac{d+1}{2q}+\frac \alpha q)}$.  
Since  $ \|f\|_{L^p(\mathbb{R}^d)} \sim 2^{j\frac{d+1}{2}(1-\frac1p)}$,  the estimate \eqref{hohoho} implies
 \[  2^{j(\frac{d+1}2 -\frac{d+1}{2q}+\frac \alpha q)} \lesssim  2^{sj} 2^{j\frac{d+1}{2}(1-\frac1p)} \] 
 for all $j\ge 1$. Thus, letting $j\to \infty$ gives   $s\geq  s_2(p,q)$.

\noindent \textit{Lower bound $s\geq s_3(p,q):= \frac{d}{p}-\frac{1-\alpha}{q}-\frac{d-1}{2}$.}
Choose again a set  $ \mathcal{E} $ as in Lemma \ref{lem:existA}. Let us consider $f$ given by 
\[  \widehat{f}(\xi)=\beta_1(|\xi|/2^{j}).\]
It is clear that $ \|f\|_{L^p(\mathbb{R}^d)} \sim 2^{jd/p'}$.  

Recall  \eqref{hhh} and \eqref{hh}, which are used to show  the  first lower bound  $s\geq  s_1(p,q)$. Translation $t\to t+1$ 
yields 
$ |e^{it\sqrt{-\Delta}}f(x)|\gtrsim 2^{j(d+1)/2} $ for  $|x| \in  [t-c_02^{-j},t]$ and $ t\in [1,2]$.  
Thus, $\|e^{it\sqrt{-\Delta}}f\|_{L^q(\mathbb{R}^d)}\gtrsim 2^{j(\frac{d+1}2 -\frac{1}{q})}$ for $ t\in [1,2]$.  Moreover, since $\#\mathcal E \sim 2^{\alpha j}$, the estimate \eqref{smooth-pq}  implies 
\[ 2^{j(\frac{d+1}2 -\frac{1}{q}+\frac \alpha q)} \lesssim 2^{js}2^{jd/p'},\] which leads to the lower bound $ s\geq s_3(p,q)$.
 
\noindent \textit{$\epsilon$-loss for  the marginal point  $ (p,q)=(\frac{2(\alpha+d-1)}{d-1},\frac{2(\alpha+d-1)}{d-1})$.} 
Let $p=q=\frac{2(\alpha+d-1)}{d-1} $ and $ \alpha>0 $. 
We again take the $2^{-j}$-separated set $\mathcal{E} $, which is  constructed in the proof of Lemma~\ref{lem:existA}. Choosing the function $f$ given by \eqref{hhh}, we have
\eqref{hh} for $t\in \mathcal{E}$ and $|x| \in J_t$. Since $\mathcal{E} $ is $2^{-j}$ separated and $ q=\frac{2(\alpha+d-1)}{d-1} $,  it follows that
\begin{align*}
\|e^{it\sqrt{-\Delta}}f\|^q_{L^q(\mathbb{R}^d\times \mathcal{E})}&\gtrsim 2^{j\frac{d+1}{2}q} \sum_{t\in\mathcal{E}}\int_{J_t}r^{-\frac{d-1}{2}q}r^{d-1}\,dr\\
&\simeq 2^{j\frac{d+1}{2}q}2^{-j} \sum_{t\in\mathcal{E}}(t-1)^{-(d-1)(\frac{q}{2}-1)}\\
& \simeq2^{j(\frac{(d+1)q}{2}-1)}\sum_{t\in\mathcal{E}}|t-1|^{-\alpha}.
\end{align*}
Recall that $\mathcal{E}$ can be written as $ \mathcal{E} =\bigcup_{l=0}^k\mathcal{E}_l$ (see Remark \ref{rem:description}) and note that $ t-1\sim \mu^{l} $ when $ t\in \mathcal{E}_l $. Thus,  we have
\[\sum_{t\in\mathcal{E}}|t-1|^{-\alpha}=\sum_{l=0}^{k}\sum_{t\in\mathcal{E}_l}|t-1|^{-\alpha}\sim \sum_{l=0}^k\#\mathcal{E}_l\mu^{-l\alpha}\sim \sum_{l=0}^k2^{k-l}2^l\sim k2^k.\]
Recalling $ k=j\alpha+O(1) $, we see that \eqref{hohoho} implies 
\[2^{j(\frac{d+1}{2}-\frac{1}{q}+\frac{\alpha}{q})}j^{1/q}\lesssim  2^{j(s+\frac{d+1}{2}-\frac{1}{p})}. 
\]
Letting $ j\to\infty $ shows that $s>{\alpha}/{q}$.
\end{proof}

\subsection*{Sharpness of the exponent $s \ge s_c(p,q)$ for \eqref{smooth-pq}}
We can now show that the regularity exponent $ s \ge s_c(p,q) $ if the local smoothing estimate \eqref{smooth-pq} holds for every $ 2^{-j} $-separated subset $ E(2^{-j}) \subset E $, whenever $ E \subset [1,2] $ has bounded $ \alpha $-Assouad characteristic. This follows from the arguments used in the proof of Proposition~\ref{prop:neces}.

In particular, let $ I_{k,1}^\mu, \ldots, I_{k,2^k}^\mu $ be the $ 2^k $ intervals of length $ \mu^k $ appearing at the $ k $-th stage of the construction used in the proof of  Lemma~\ref{lem:existA}. We consider the set
\[
E = \bigcap_{k=0}^\infty \bigcup_{\ell=1}^{2^k} I_{k,\ell}^\mu.
\]
It is easy to see that $ E \subset [1,2] $ has bounded $ \alpha $-Assouad characteristic. Recall that the set $ \mathcal{E} $ consists of the right endpoints of the $ 2^k $ intervals $ I_{k,1}^\mu, \ldots, I_{k,2^k}^\mu $. Clearly, $ \mathcal{E} $ is a $ 2^{-j} $-separated subset of $ E $. Therefore, by the same reasoning as in the proof of Proposition~\ref{prop:neces}, we conclude that $ s \ge s_c(p,q) $.

 \section*{Acknowledgement} 
This work was supported by the National Research Foundation of Korea (RS-2024-00342160; Lee and Zhao); 
the Basque Government  and Spanish MICIN and MICIU (BERC programme and Ikerbasque,  CEX2021-001142-S, PID2023-146646NB-I00, CNS2023-143893; Roncal); 
and the China Scholarship Council (Zhang).

\bibliographystyle{plain}

\end{document}